\documentclass[10pt]{amsart}
\usepackage[pagebackref]{hyperref}
\usepackage{amscd}

\renewcommand{\Re}{\operatorname{Re}}
\renewcommand{\Im}{\operatorname{Im}}
\renewcommand{\d}{\mathrm{d}}

\newcommand{\I}{\mathrm{I}}
\newcommand{\e}{\mathrm{e}}
\newcommand{\iC}{\mathrm{i}}

\newcommand{\ts}{\textstyle }

\newcommand{\ph}{\phantom }

\newcommand{\cH}{{\mathcal H}}
\newcommand{\cK}{{\mathcal K}}
\newcommand{\cM}{{\mathcal M}}
\newcommand{\cR}{{\mathcal R}}

\newcommand{\bbR}{{\mathbb R}}
\newcommand{\R}[1]{\ensuremath{\mathbb R^{\,#1}}} % to fix supscript spacing

\newcommand{\bbC}{{\mathbb C}}
\newcommand{\C}[1]{\ensuremath{\mathbb C^{\,#1}}} % to fix supscript spacing
\newcommand{\bbP}{{\mathbb P}}
\newcommand{\tbbP}{\tilde{\mathbb P}}
\newcommand{\bbS}{{\mathbb S}}

\newcommand{\ov}{\overline}
\newcommand{\ot}{\otimes}

\newcommand{\om}{\omega}
\newcommand{\tht}{\theta}
\newcommand{\eb}{\mathbf{e}}
\newcommand{\fb}{\mathbf{f}}

\newcommand{\eug}{\operatorname{\mathfrak{g}}}
\newcommand{\eugl}{\operatorname{\mathfrak{gl}}}
\newcommand{\euh}{\operatorname{\mathfrak{h}}}
\newcommand{\eut}{\operatorname{\mathfrak t}}

\newcommand{\GL}{\operatorname{GL}}
\newcommand{\SL}{\operatorname{SL}}
\newcommand{\SO}{\operatorname{SO}}
\newcommand{\Un}{\operatorname{U}}
\newcommand{\Or}{\operatorname{O}}
\newcommand{\Gr}{\operatorname{Gr}}

\DeclareMathOperator{\tr}{tr}

\newcommand{\w}{{\mathchoice{\,{\scriptstyle\wedge}\,}{{\scriptstyle\wedge}}
      {{\scriptscriptstyle\wedge}}{{\scriptscriptstyle\wedge}}}}
\newcommand{\lhk}{\mathbin{\hbox{\vrule height1.4pt width4pt depth-1pt 
             \vrule height4pt width0.4pt depth-1pt}}}

\numberwithin{equation}{subsection}

\newtheorem{theorem}{Theorem}

\newtheorem{proposition}{Proposition}
\newtheorem{corollary}{Corollary}

\theoremstyle{remark}
\newtheorem{definition}{Definition}
\newtheorem{remark}{Remark}
\newtheorem{example}{Example}

\begin{document}

\author[R. Bryant]{Robert L. Bryant}
\address{Duke University Mathematics Department\\
         P.O. Box 90320\\
         Durham, NC 27708-0320}
\email{\href{mailto:bryant@math.duke.edu}{bryant@math.duke.edu}}
\urladdr{\href{http://www.math.duke.edu/~bryant}%
         {http://www.math.duke.edu/\lower3pt\hbox{\symbol{'176}}bryant}}

\dedicatory{This article is dedicated to Shiing-Shen Chern, 
whose beautiful works on Finsler geometry have inspired
so much progress in the subject.}

\title[]
      {Some remarks on Finsler manifolds  \\
          with constant flag curvature }

\date{July 31, 2001}

\begin{abstract}
This article is an exposition of four loosely related
remarks on the geometry of Finsler manifolds with 
constant positive flag curvature.  

The first remark
is that there is a canonical K\"ahler structure on the
space of geodesics of such a manifold.  

The second remark is that there is a natural way to construct a 
(not necessarily complete) Finsler $n$-manifold of constant 
positive flag curvature out of a hypersurface in suitably general
position in~$\bbC\bbP^n$.

The third remark is that there is a description of the Finsler
metrics of constant curvature on~$S^2$ in terms of a Riemannian
metric and $1$-form on the space of its geodesics.  In particular,
this allows one to use any (Riemannian) Zoll metric of positive 
Gauss curvature on~$S^2$ to construct a global Finsler metric 
of constant positive curvature on~$S^2$.

The fourth remark concerns the
generality of the space of (local) Finsler metrics 
of constant positive flag curvature in dimension~$n{+}1>2$.  
It is shown that such metrics depend on $n(n{+}1)$
arbitrary functions of $n{+}1$ variables and that 
such metrics naturally correspond to certain torsion-free
$S^1{\cdot}\GL(n,\bbR)$-structures on $2n$-manifolds.
As a by-product, it is found that these groups do occur
as the holonomy of torsion-free affine connections 
in dimension~$2n$, a hitherto unsuspected phenomenon.
\end{abstract}



\thanks{
Thanks are due to Duke University 
for its support via a research grant 
and to the National Science Foundation 
for its support via DMS-9870164.  This article
was finished during a July 2001 visit to IMPA;
I thank them for their hospitality.
\hfill\break
\hspace*{\parindent} 
This is Version~$1.1$. 
The most recent version can be found at arXiv:math.DG/0107228 .
}

\maketitle

\setcounter{tocdepth}{2}
\tableofcontents

\section[Introduction]{Introduction}\label{sec: intro}

The purpose of this article is to explain some new results
in the theory of Finsler manifolds with constant flag 
curvature, particularly constant positive flag curvature.

For general background in the subject, the reader can consult
\cite{BCS, Ma2, Ru} and for articles dealing specifically 
with the case of constant flag curvature, the reader may 
consult~\cite{Ak, Fu3, Ma1, Sh}.

\subsection{The main results}\label{ssec: intro main results}
Though the discussion in this article will hold for a wider
notion of Finsler structure than is usually considered, the
statements made in this introduction will be focussed on the
case of a classical (though not necessarily reversible) Finsler
structure on a manifold.  

Suppose that~$M$ is an $(n{+}1)$-manifold endowed with 
a Finsler structure, regarded as being specified by its 
unit tangent bundle~$\Sigma\subset TM$ (often referred to 
as the \emph{tangent indicatrix}).  Suppose further 
that~$M$ is geodesically simple, i.e., that the quotient~$Q$
of~$\Sigma$ by the geodesic flow can be given the structure
of a smooth $2n$-manifold in such a way that the quotient 
map~$q:\Sigma\to Q$ is a smooth submersion.%
\footnote{
Of course, this can always be arranged locally by restricting
attention to a geodesically convex neighborhood in~$M$.}

As is well-known in symplectic geometry, the space~$Q$, 
which can be thought of as the space of oriented geodesics 
of the Finsler structure, inherits a canonical symplectic structure.

According to Theorem~\ref{thm: Kahler str on Q}, when the
Finsler structure has constant positive flag curvature, $Q$
also inherits a natural Riemannian metric with respect to which
the symplectic form is parallel.  In other words, $Q$ is naturally
a K\"ahler manifold.

It turns out that~$Q$ has a yet finer structure. For
each~$x\in M$, the set~$Q_x\subset Q$ consisting of the geodesics
that pass through~$x$ is a totally real submanifold of~$Q$. 
For a fixed geodesic~$q\in Q$, the set of manifolds~$Q_x$ 
as~$x\in M$ varies on~$q$ defines a $1$-parameter family of totally
real submanifolds of~$Q$ passing through~$q$.  In the case
that the Finsler structure has constant flag curvature~$1$,
the totally real tangent planes~$T_qQ_x\subset T_qQ$ 
as $x$ varies over~$q$ turn out
to differ by multiplication by complex numbers of
the form~$\e^{\iC\theta}$, i.e., there is a canonical circle
of totally real $n$-planes passing through each point of~$Q$.
This defines a canonical~$S^1{\cdot}\Or(n)$-structure on~$Q$.
This $S^1{\cdot}\Or(n)$-structure is \emph{not} torsion-free 
except in the trivial case where~$M$ is a Riemannian manifold 
of constant positive sectional curvature.

However, as is shown in~\S\ref{sssec: S1dotGLnR str}, 
this $S^1{\cdot}\Or(n)$-structure on~$Q$ underlies
a canonical $S^1{\cdot}\GL(n,\bbR)$-structure that \emph{is}
torsion-free.  This is surprising,
since, for~$n>2$, the group~$S^1{\cdot}\GL(n,\bbR)\subset\GL(2n,\bbR)$
was not previously recognized to be possible as holonomy of
a torsion-free connection on a $2n$-manifold.  Nevertheless, 
as Theorem~\ref{thm: structure generality} shows, these groups 
are indeed realizable as holonomy groups in this way.

In fact, it turns out (\S\ref{sec: generality}) that there is a very close 
connection between torsion-free $S^1{\cdot}\GL(n,\bbR)$-structures 
on $2n$-manifolds
and Finsler structures with constant flag curvature~$1$.  When~$n>2$,
a torsion-free $S^1{\cdot}\GL(n,\bbR)$-structure on a $2n$-manifold~$Q$
that satisfies a mild positivity condition on its curvature arises 
from a canonical (generalized) Finsler structure 
of constant flag curvature~$1$ on an $(n{+}1)$-manifold~$M$.
When~$n=2$, one must impose a further condition on the torsion-free
structure, that of \emph{integrability}, but the local generality
of the integrable, torsion-free $S^1{\cdot}\GL(2,\bbR)$-structures
is also easily analyzable from this standpoint.
Thus, the construction is reversible, so that 
Theorem~\ref{thm: structure generality}
gives a method of describing the local generality of (generalized)
Finsler structures of constant flag curvature. 

The other main results deal with either special dimensions
or more special Finsler structures:

First of all, an old result of Funk~\cite{Fu3} describes the local 
Finsler metrics on the plane that have constant positive curvature
and are rectilinear (i.e., the geodesic paths are straight lines)
in terms of a holomorphic function of one variable.  In~\cite{Br3},
this construction was given a projectively invariant description
in terms of certain holomorphic curves without real points in~$\bbC\bbP^2$.
This turns out to generalize in a natural way to higher dimensions:
A Finsler metric on a domain in~$\bbR^{n+1}$ with constant flag curvature 
whose geodesics are straight lines gives rise to a holomorphic 
hypersurface~$Q\subset\bbC\bbP^{n+1}$ satisfying certain open conditions
and, conversely, such a hypersurface determines a (generalized) Finsler
structure on a domain in~$\bbR^{n+1}$ in a projectively natural way.
For a precise statement, see Theorem~\ref{thm: convex gives CFC}.

This result is used to derive two further results:  First, the global 
Finsler metrics on~$\bbR\bbP^{n+1}$ with constant flag curvature~$1$
and rectilinear geodesics are determined 
(Example~\ref{ex: non-real quadrics}).  It turns out that these correspond
naturally to the hyperquadrics in~$\bbC\bbP^{n+1}$ that have no real
points.  Thus, up to isomorphism, these consist of an $(n{+}1)$-parameter
family of distinct global examples.  Second, it is shown that for any 
closed, real analytic hypersurface~$S\subset T_0\bbR^{n+1}\simeq\bbR^{n+1}$
that is strictly convex towards the origin, there exists a Finsler metric
with constant flag curvature~$1$ and rectilinear geodesics 
on a neighborhood~$U$ of~$0\in\bbR^{n+1}$ that has~$S$ as its space of
unit tangent vectors at~$0\in\bbR^{n+1}$.  

Finally, in~\S\ref{ssec: dim 2} 
the description in~\cite{Br1} of Finsler metrics on~$S^2$ 
of constant positive curvature~$1$ in terms of a Riemannian 
metric~$\d\sigma^2$ and a (`magnetic') $1$-form~$\beta$ on~$Q\simeq S^2$,
the space of oriented geodesics of the Finsler structure, is recalled
and then combined with Guillemin's classic result on the existence
of Zoll metrics on the $2$-sphere to prove the existence 
of a large family of global Finsler metrics on~$S^2$ 
with constant positive curvature~$1$.  This is still far from a complete
description, of course, but it gives an indication that this family
is much larger than previously believed.

\subsection{Dedication}\label{ssec: dedication}
This article is dedicated to Professor Shiing-Shen Chern, 
whose research in and tireless efforts to promote the study 
of Finsler geometry for nearly 60 years have inspired much 
of the progress in the subject.  The research in this article
would not have been possible without his encouragement and
interest.

The main results in this article (aside from those 
of~\S\ref{sec: generality}) were announced at the 1998 
Geometry Festival at Stony Brook, NY, but I had been unable 
to set aside time to write the article until this opportunity 
came along. Thus, I would like to thank the editorial board 
of the Houston Journal of Mathematics for inviting me to write 
this article for a volume dedicated to Professor Chern.  

\section[The Structure Equations]{The Structure Equations}
\label{sec: str eqs}
This first section is mainly to fix notation and to 
remind the reader of some basic facts about Finsler
geometry that will be used in this article.  It will
also be necessary to generalize the notion of Finsler
structure slightly since some of the constructions that
will be made have to first be done in this slightly
more general context.

\subsection{Generalized Finsler structures}
\label{ssec: gen Finsler strs}
Let~$M$ be a manifold of dimension~$n{+}1$.  Classically,
a Finsler structure on~$M$ is a non-negative function~$F:TM\to\bbR$ 
that is smooth and positive away from the zero section of~$TM$, 
homogeneous of degree~$1$ (i.e., $F(\lambda v) = |\lambda|\,F(v)$
for all~$v\in TM$ and~$\lambda\in\bbR$), and strictly convex
on each tangent space~$T_xM$ for~$x\in M$.  For background,
the reader is referred to~\cite{BCS}.  

The function~$F$ determines and is determined by the set
\begin{equation}\label{eq: F indicatrix}
\Sigma_F = \{ v\in TM\mid F(v)=1\},
\end{equation}
which is known as the \emph{tangent indicatrix} or \emph{unit
tangent bundle} of~$F$.  For each~$x\in M$, the 
intersection~$\Sigma_F(x) = \Sigma_F\cap T_xM$ is a smooth,
compact hypersurface in the vector space~$T_xM$ that is
transverse to the radial vector field on~$TM$ and is 
strictly convex towards the origin.

\begin{definition}\label{def: gen Finsler str}
A \emph{generalized Finsler structure} on a manifold~$M^{n+1}$
is a pair~$(\Sigma,\iota)$ where~$\Sigma$ is a connected, smooth manifold
of dimension~$2n{+}1$ together with a
radially transverse immersion~$\iota:\Sigma\to TM$
with the following two properties:
\begin{enumerate}
\item The composition $\pi\circ\iota:\Sigma\to M$ is a submersion
      with connected fibers.
\item Setting~$\Sigma_x= \iota^{-1}(T_xM)$ for each~$x\in M$,
the mapping~$\iota_x:\Sigma_x\to T_xM$ immerses~$\Sigma_x$ as
a hypersurface in~$T_xM$ that is locally strictly convex towards
the origin~$0_x$.
\end{enumerate}
\end{definition}

\begin{remark}[Equivalence]\label{rem: equivalence}
Two generalized Finsler structures, say $(\Sigma_1,\iota_1)$ on~$M_1$ 
and~$(\Sigma_2,\iota_2)$ on~$M_2$, will be said to be \emph{isometric} 
if there is a diffeomorphism~$\psi:\Sigma_1\to\Sigma_2$
and a diffeomorphism~$\phi:M_1\to M_2$
so that $\phi'\circ\iota_1 = \iota_2\circ\phi$.  

The reader might prefer to regard a generalized Finsler 
structure as an isometry class of generalized Finsler structures 
as they were defined in Definition~\ref{def: gen Finsler str}.  While
this is natural, it can be cumbersome, so this course has not been
adopted.
\end{remark}

Of course, the canonical inclusion into~$TM$ of the tangent indicatrix 
of a classical Finsler structure on~$M$ is a generalized Finsler
structure.  Obviously, this is not the only kind
of example.  For instance, there is no requirement that any
of the~$\Sigma_x$ be compact, or that~$\iota$ be an embedding.

The reasons for considering generalized Finsler structures is two-fold.
First, all of the classical constructions of canonical connections,
bundles, and curvature will work just as well for generalized Finsler
structures as for Finsler structures with no increase in difficulty.  
Second, as will be seen, imposing
differential equations (such as curvature constraints) on Finsler
structures often leads to problems where the best strategy is to
first solve the problem in the more general class of generalized
Finsler structures and then look among the solutions for Finsler
structures in the classical sense. 

\subsection{The structure bundle}\label{ssec: str bndl for FS}
Let~$(\Sigma,\iota)$ be a generalized Finsler structure on a
manifold~$M^{n+1}$. Following Cartan~\cite{Ca1,Ca2} 
and Chern~\cite{Ch1,Ch2}, one can define a canonical~$\Or(n)$-structure
with connection on~$\Sigma$.  This section will review 
their constructions via the method of equivalence 
and establish the notation to be used throughout this article.

\subsubsection{The Hilbert form}\label{sssec: hilbert form}
One constructs a contact form~$\om_0$ on~$\Sigma$
as follows:  For each~$u\in\Sigma$, the vector~$\iota(u)$ lies in~$T_xM$
where~$x = \pi\bigl(\iota(u)\bigr)$.  Moreover, the 
image~$(\pi{\circ}\iota)'(T_u\Sigma_x) \subset T_xM$ is a hyperplane
not containing~$\iota(u)$.  Consequently, there exists a unique linear
form~$u^*\in T^*_xM$ whose kernel is~$(\pi{\circ}\iota)'(T_u\Sigma_x)$
and so that~$u^*\bigl(\iota(u)\bigr)=1$.  Define 
the $1$-form~$\om_0$ on~$\Sigma$ so that
\begin{equation}\label{eq: Hilbert defined}
(\om_0)_u = (\pi{\circ}\iota)^*(u^*).
\end{equation}

The assumption that~$\iota:\Sigma\to TM$ is radially transverse
(i.e., transverse to the orbits of scalar multiplication on~$TM$) 
implies that~$\om_0$ is a contact
form, i.e., that~$\om_0\w(\d\om_0)^n\not=0$.  This form is
known in the calculus of variations as the \emph{Hilbert form}.

\subsubsection{The Reeb field}\label{sssec: Reeb field}
Since~$\om_0$ is a contact form, 
there exists a unique vector field~$E$ on~$\Sigma$ that satisfies
$\om_0(E) = 1$ and~$E\lhk(\d\om_0) = 0$.  This vector field is 
known as the \emph{Reeb vector field}.  Its flow on~$\Sigma$
is simply the geodesic flow when~$(\Sigma,\iota)$ is an actual
Finsler structure, so it will be referred to as the geodesic
flow or the Reeb flow in this more general context.

In particular, a $\Sigma$-geodesic will be a smooth
curve~$\gamma:(a,b)\to M$ such that~$\gamma':(a,b)\to TM$
lifts back to~$\Sigma$ as an integral curve of~$E$.

The generalized Finsler structure~$\Sigma$ will be said to be
\emph{geodesically complete} if~$E$ is complete, i.e., if the
flow of~$E$ is globally defined.

\subsubsection{The Legendrian foliation}\label{sssec: Leg fol}
The foliation~$\cM$ whose leaves 
are the fibers~$\Sigma_x$ for~$x\in M$ is $\om_0$-Legendrian.
As a consequence, each point~$u\in\Sigma$ has a neighborhood, say~$U$, 
on which there exist $n$ $1$-forms~$\om_1,\dots,\om_n$
with the properties that 
\begin{enumerate}
\item $\om_0\w\dots\w\om_n\not=0$,
\item each of the~$\om_i$ vanishes when pulled back to any~$\Sigma_x$,
\item $\om_i(E)=0$ for~$1\le i\le n$, and
\item $\d\om_0\w\om_1\w\dots\w\om_n=0$.
\end{enumerate}

It follows that there exist~$1$-forms~$\tht^1,\dots,\tht^n$ on~$U$ so that
\begin{equation}\label{eq: d omega_0 tight}
\d\om_0 = {} -  \tht^1\w\om_1 - \dots - \tht^n\w\om_n\,.
\end{equation}
The forms~$\om_0,\om_1,\dots,\om_n,\tht^1,\dots,\tht^n$ 
are linearly independent on~$U$.

The conditions imposed on the $n$ $1$-forms~$\om_1,\dots,\om_n$ 
so far determine them up to a change of basis (with variable coefficients).
If one were to make a different choice subject to the same conditions,
one would have $1$-forms
\begin{equation}\label{eq: star om}
{}^*\!\om_i = A^j_i\,\om_j
\end{equation}
for some invertible matrix of functions~$A = (A^i_j)$ defined on~$U$.  
Any~${}^*\tht^1,\dots,{}^*\tht^n$ for which the
equation~$\d\om_0 = -{}^*\!\tht^i\w {}^*\!\om_i$ 
holds must then be of the form
\begin{equation}\label{eq: star tht}
{}^*\!\tht^i = B^i_j\,(\tht^j + S^{jk}\,\om_k)
\end{equation}
where~$(B^i_j) = B = {}^t\!A^{-1}$ and~$S^{ij} = S^{ji}$ 
are functions on~$U$. 

In the language of~$G$-structures, 
the local coframings~$(\om_0,\om_i,\tht^i)$ that satisfy
the above conditions are the local sections of a~$G_1$-bundle
over~$\Sigma$ where~$G_1\subset \GL(2n{+}1,\bbR)$ is the
group of matrices of the form
\begin{equation}\label{eq: G_1 matrices}
\left\{\ \left[\begin{matrix}
1&0&0\\ 0 & A & 0\\ 0 & {}^t\!A^{-1}S & {}^t\!A^{-1}
\end{matrix}\right]
\quad\vrule\quad 
A\in\GL(n,\bbR),\ \ S={}^t\!S\in M_n(\bbR)\ 
\right\}\,.
\end{equation}
Such coframings are said to be \emph{$1$-adapted}.

\subsubsection{Convexity}\label{sssec: convexity}
Since the system spanned by~$\om_0,\om_1,\dots,\om_n$ is Frobenius,
there must be functions~$H_{ij}$ on~$U$ so that
\begin{equation}\label{eq: d om_i}
\d\om_i \equiv H_{ij}\,\tht^j\w\om_0 \mod \om_1,\dots,\om_n\,.
\end{equation}

Using~\eqref{eq: d omega_0 tight} to expand the 
identity~$\d(\d\om_0)=0$, reducing modulo~$\om_1,\dots,\om_n$, 
and then using~\eqref{eq: d om_i} shows that~$H_{ij}=H_{ji}$. 

The geometric significance of~$H$ is that the 
quantity~$H_{ij}\,\tht^i{\circ}\tht^j$ pulls back
to each~$\Sigma_x$ to be the centro-affine invariant metric
induced on it by its radially transverse immersion into
the vector space~$T_xM$. 

In particular, the strict local convexity hypothesis
implies that the symmetric matrix~$H=(H_{ij})$ 
is positive definite everywhere on~$U$.

Moreover, relative to a coframing~$\bigl(\om_0,{}^*\!\om_i,{}^*\!\tht^i\bigr)$
on~$U$ given by~\eqref{eq: star om} and~\eqref{eq: star tht}, 
the corresponding symmetric matrix~${}^*\!H$ satisfies
${}^*\!H = A\,H\,{}^t\!A$.  

Thus, the coframing~$\om_0,\om_1,\dots,\om_n,\tht^1,\dots,\tht^n$
can be chosen so that it satisfies~$H=\I_n$, i.e.,
so that
\begin{equation}\label{eq: d om_i sharpened}
\d\om_i \equiv \,\tht^j\w\om_0 \mod \om_1,\dots,\om_n\,.
\end{equation}

Henceforth, assume that~\eqref{eq: d omega_0 tight} 
and~\eqref{eq: d om_i sharpened} hold.  
In the language of~$G$-structures, 
the local coframings~$(\om_0,\om_i,\tht^i)$ 
that satisfy these conditions are the local sections of a~$G_2$-bundle
over~$\Sigma$ where~$G_2\subset \GL(2n{+}1,\bbR)$ is the
group of matrices of the form
\begin{equation}\label{eq: G_2 matrices}
\left\{\ \left[\begin{matrix}
1&0&0\\ 0 & A & 0\\ 0 & AS & A
\end{matrix}\right]
\quad\vrule\quad 
A\in\Or(n),\ \ S={}^t\!S\in M_n(\bbR)\ 
\right\}\,.
\end{equation}
Such coframings are said to be \emph{$2$-adapted}.

\subsubsection{The quadratic form~$\gamma$}\label{sssec: gamma defined}
The forms~$\om_1,\dots,\om_n$ now are determined 
up to an orthogonal change of basis.
Thus, the quadratic form
\begin{equation}\label{eq: quadratic form}
\gamma = {\om_1}^2 + \dots + {\om_n}^2
\end{equation}
is globally well-defined on~$\Sigma$.  
Here is the geometric meaning of~$\gamma$:
 
For any~$u\in\Sigma$ with basepoint~$x = \pi\bigl(\iota(u)\bigr)$
in~$M$, there is a unique positive definite quadratic 
form~$g_u$ on~$T_xM$ with the property that~$\iota(u)$ is a unit vector
for~$g_u$ and that the unit sphere of~$g_u$ in~$T_xM$ osculates
to second order to~$\iota(\Sigma_x)$ at~$\iota(u)$.  This family
of quadratic forms can be shown to satisfy
\begin{equation}\label{eq: gamma characterized}
(\pi{\circ}\iota)^*(g_u) 
  = \left({\om_0}^2 + \gamma\right)\vrule_u
  = \left({\om_0}^2 + {\om_1}^2 + \dots + {\om_n}^2\right)\vrule_u\,.
\end{equation}

Since the ambiguity in the choice of the~$\om_i$
lies in the orthogonal group, there is no longer any reason to
preserve a distinction between upper and lower indices.  
\emph{Henceforth, all indices will be written as subscripts.}
In particular,~$\tht^i$ will now be written as~$\tht_i$.

\subsubsection{Further normalizations}\label{sssec: further normalizations}
In view of~\eqref{eq: d om_i sharpened}, there must exist 
$1$-forms~$\om_{ij}$ on~$U$ so that
\begin{equation}\label{eq: d om_i sharpened more}
\d\om_i = \,\tht_j\w\om_0 - \om_{ij}\w\om_j\,.
\end{equation}
The relations~\eqref{eq: d om_i sharpened more} do not
determine the~$\om_{ij}$ uniquely.  Evidently, one can
keep the same relations while replacing
each~$\om_{ij}$ by ${}^*\om_{ij}=\om_{ij}+P_{ijk}\,\om_k$
where~$P_{ijk}=P_{ikj}$ are arbitrary functions on~$U$.

Write~$\om_{ij}=\tht_{ij}+\sigma_{ij}$ 
where~$\tht_{ij}=-\tht_{ji}$
and~$\sigma_{ij}=\sigma_{ji}$.  Expand~$\sigma_{ij}$ in
the coframing as follows:
\begin{equation}\label{eq: sigma_ij crude}
\sigma_{ij} = S_{ij}\,\om_0 + B_{ijk}\,\om_k+ I_{ijk}\,\tht_k\,,
\end{equation}
where~$S_{ij}=S_{ji}$, $B_{ijk}=B_{jik}$, and~$I_{ijk}=I_{jik}$
are functions on~$U$.    

The ambiguities in the choices so far can be exploited 
to eliminate the quantities~$S_{ij}$ and~$B_{ijk}$ by the
following normalizations:

\emph{A Lagrangian splitting.} First, note that by replacing~$\tht_i$ 
by~${}^*\tht_i = \tht_i - S_{ij}\,\om_j$ 
and~$\om_{ij}$ by~${}^*\om_{ij} = \om_{ij}-S_{ij}\,\om_0$, 
one preserves the formulae~\eqref{eq: d omega_0 tight}
and~\eqref{eq: d om_i sharpened}, but the~$S_{ij}$ are replaced
by~${}^*\!S_{ij}=0$.  

\emph{Thus, it will be assumed from now on that~$S_{ij}=0$.}

\begin{remark}[Geometric $S$ elimination]\label{rem: S elimination}
This normalization has the following intrinsic description,
which, in slightly different form, can
essentially be found in the work of Foulon~\cite{Fo}:

Since~$\om_0(E)=1$ and $\om_i(E)=\tht_i(E)=0$, the
Lie derivative~$\dot\gamma$ of~$\gamma$ with respect to~$E$ 
can be computed in the form
\begin{equation}\label{eq: Lie derivative gamma}
\begin{split}
\dot\gamma = 2\,\om_i\circ(E\lhk \d\om_i) 
&= 2\,\om_i\circ\bigl(-\tht_i-(\tht_{ij}(E)+\sigma_{ij}(E))\,\om_j\bigr)\\
&= -2\,\om_i\circ(\tht_i-S_{ij}\,\om_j)
\end{split}
\end{equation}
(remember that~$\tht_{ij}(E)=-\tht_{ji}(E)$).  

It follows that the $(n{+}1)$-plane field on~$U$ 
defined by~$\tht_i-S_{ij}\,\om_j=0$ is the unique one that is
transverse to the fibers of~$\pi\circ\iota$ and is
both Lagrangian with respect to~$\d\om_0$ and null with respect
to the quadratic form~$\gamma'$.  Consequently, this plane field
is globally defined on~$\Sigma$.
\end{remark}

In the language of~$G$-structures, 
the $2$-adapted coframings~$(\om_0,\om_i,\tht^i)$ 
that satisfy equations of the form
\begin{equation}
\d\om_i = \tht_i\w\om_0 - \tht_{ij}\w\om_j 
- B_{ijk}\,\om_k\w\om_j - I_{ijk}\,\tht_k\w\om_j\,,
\end{equation}
where~$\tht_{ij}=-\tht_{ji}$ and~$B_{ijk}=B_{jik}$ 
while~$I_{ijk}=I_{jik}$,
are the local sections of a~$G_3$-bundle
over~$\Sigma$ where~$G_3\subset \GL(2n{+}1,\bbR)$ is the
group of matrices of the form
\begin{equation}\label{eq: G_3 matrices}
\left\{\ \left[\begin{matrix}
1&0&0\\ 0 & A & 0\\ 0 & 0 & A
\end{matrix}\right]
\quad\vrule\quad 
A\in\Or(n)\ 
\right\}\,.
\end{equation}
Such coframings are said to be \emph{$3$-adapted}.

\emph{A connection adaptation.} By making a 
replacement~$\tht_{ij}\longmapsto{}^*\tht_{ij}=\tht_{ij}+P_{ijk}\,\om_k$
where~$P_{ijk}=-P_{jik}$, one can arrange~${}^*\!B_{ijk}=0$ 
and this uniquely determines the~$P_{ijk}$.%
\footnote{ This is precisely
the algebraic lemma that is used to prove the Fundamental Lemma
of Riemannian geometry.}
Thus, it will be assumed from now on that~$B_{ijk}=0$.  
This normalization determines the~$\tht_{ij}$ uniquely.

\subsubsection{An $\Or(n)$-structure}\label{sssec: O(n) structure}
The analysis so far has produced 
a coframing
\begin{equation}\label{eq: coframing list}
(\om_0,\om_1,\dots,\om_n,\tht_1,\dots,\tht_n)
\end{equation}
on~$U$ that satisfies the equations
\begin{equation}\label{eq: structure summary}
\begin{split}
\d\om_0 &=  \ph{\tht_j\w\om_0 -\tht_{ij}}
      \hbox to 0pt{ \hss ${} - \tht_j$}\w\om_j\,\\
\d\om_i &=  \tht_i\w\om_0 -\tht_{ij}\w\om_j 
               - I_{ijk}\,\tht_{k}\w\om_j
\end{split}
\end{equation}
for some functions~$I_{ijk}=I_{jik}$ and $1$-forms~$\tht_{ij}=-\tht_{ji}$.  
Moreover, coframings
satisfying these equations are unique up to an orthogonal change 
of coframing of the form~${}^*\!\om_i = A_{ij}\,\om_j$ 
and~${}^*\!\tht_i = A_{ij}\,\tht_j$ where~$A = (A_{ij})$ is a 
smooth mapping of~$U$ into~$\Or(n)$.

Thus, such coframings are the local sections 
of an $\Or(n)$-structure~$u:F\to\Sigma$, 
where~$\Or(n)$ is embedded into $\GL(2n{+}1,\bbR)$ as the
subgroup~$G_3$ of matrices of the form~\eqref{eq: G_3 matrices}.
 
The~$\tht_{ij}=-\tht_{ji}$
are simply connection forms for this~$\Or(n)$-structure
relative to the given coframing.%
\footnote{ This is essentially Cartan's connection, but, as was observed 
explicitly by Chern~\cite{Ch2}, there are other natural connections 
one could conceivably attach to~$F$. For example, one could
take the Levi-Civita connection of the Riemannian 
metric~$\d s^2 = {\om_0}^2+{\om_i}^2+{\tht_i}^2$ on~$\Sigma$
and project it to~$F$ (which is a subbundle of the~$\d s^2$-orthonormal
frame bundle) in the usual way.  Indeed, several different 
connections have been attached to Finsler geometry in the literature; 
having to sort through all of them while learning the subject is 
something of a chore.  See~\cite{BCS} for an account.}

To keep the notation simple, the symbols~$\om_0,\om_i,\tht_i$ will
also be used to stand for the corresponding tautological forms on~$F$ 
while the symbols~$\tht_{ij}$ will also be used to stand for 
the connection forms on~$F$.  Context will be used to determine
whether these forms are to be understood as defined globally on~$F$ 
or locally on~$\Sigma$.  In most cases, this will make no practical
difference.   

For example, the quadratic form~$\gamma$ is defined
globally on~$\Sigma$ while the expression~${\om_1}^2+\dots+{\om_n}^2$
is defined globally on~$F$.  Logically, one should write $u^*(\gamma)
= {\om_1}^2+\dots+{\om_n}^2$ as a global equation on~$F$, but, as is
common practice in moving frame computations, one simply writes~$\gamma
= {\om_1}^2+\dots+{\om_n}^2$ and either the $u^*$ is understood or 
else the equation is meant locally on~$\Sigma$, relative to a $2$-adapted
coframing.

\subsubsection{The symmetry of~$I$}\label{sssec: I symmetric}
Differentiating the first equation of~\eqref{eq: structure summary} 
yields the relation
\begin{equation}\label{eq: first d eta}
0=\left(\d\tht_i+\tht_j\w\tht_{ji}+I_{ijk}\,\tht_j\w\tht_k\right)\w\om_i\,,
\end{equation}
which, in particular, implies
\begin{equation}\label{eq: first d eta mod om1..n}
\d\tht_i + \tht_j\w\tht_{ji} + I_{ijk}\,\tht_j\w\tht_k
\equiv 0 \mod \om_1,\dots,\om_n\,.
\end{equation}
Differentiating the second equation of~\eqref{eq: structure summary}
and reducing modulo~$\om_1,\dots,\om_n$ yields
\begin{equation}\label{eq: second d eta}
0\equiv\left(\d\tht_i+\tht_j\w\tht_{ji}+I_{ijk}\,\tht_k\w\tht_j\right)\w\om_0
\mod \om_1,\dots,\om_n\,,
\end{equation}
which, in particular, implies
\begin{equation}\label{eq: second d eta mod om0..n}
\d\tht_i + \tht_j\w\tht_{ji} - I_{ijk}\,\tht_j\w\tht_k
\equiv 0 \mod \om_0,\om_1,\dots,\om_n\,.
\end{equation}
However, comparing~\eqref{eq: first d eta mod om1..n} 
with~\eqref{eq: second d eta mod om0..n} 
yields~$I_{ijk}\,\tht_j\w\tht_k=0$, i.e.,~$I_{ijk}=I_{ikj}$.
In particular, $I_{ijk}$ is fully symmetric in its indices.

\subsubsection{The Cartan torsion}\label{sssec: Cartan torsion}
The symmetric cubic form
\begin{equation}\label{eq: Cartan torsion}
I = I_{ijk}\,\tht_i\,\tht_j\,\tht_k
\end{equation}
is well-defined globally on~$\Sigma$. Its geometric interpretation 
at a point~$u$ is that it measures the failure of the unit sphere
of~$g_u$ in~$T_xM$ to osculate to~$\iota(\Sigma_x)$ to third order
at~$\iota(u)$.  In fact, $I$ pulls back to each~$\Sigma_x$ to be
the classical centro-affine cubic form induced on~$\Sigma_x$ by
its radially transverse, strictly locally convex immersion 
into the vector space~$T_xM$.

It is a standard result of Cartan that the equation~$I\equiv0$ is
the necessary and sufficient condition that the image~$\iota(\Sigma)$
should be an open subset of the unit sphere bundle of a 
Riemannian metric~$g$ (necessarily unique)
defined on the open set~$\pi{\circ}\iota(\Sigma)\subset M$.  

In fact, note that, when~$I\equiv0$, if one writes~$\tht_i = \tht_{0i}$
and sets~$\tht_{i0}=-\tht_{0i}$ and~$\tht_{00}=0$, 
then~\eqref{eq: structure summary} can be written in the simple form
\begin{equation}\label{eq: I=0 case}
\d\om_a = {} -\tht_{ab}\w\om_b\,,\qquad\qquad 
\tht_{ab}+\tht_{ba} = 0\,,
\end{equation}
where the indices~$a$ and~$b$ lie in the range~$0,\le a,b\le n$.
It follows that the quadratic form~${\om_0}^2+{\om_1}^2+\dots+{\om_n}^2$
is simply the $\pi{\circ}\iota$-pullback of the desired metric~$g$.

\subsubsection{Realizations}\label{sssec: realizations}
Because the notational change in the Riemannian case is so suggestive,
it will be adopted here for the general case.  Thus, from now on,
$\tht_i$ will be written as~$\tht_{0i}$ while~$\tht_{i0}$
(to be used on rare occasions) will mean~$-\tht_{0i}=-\tht_i$.
With this notational change, the structure equations so far take the
form
\begin{equation}\label{eq: structure summary new}
\begin{split}
\d\om_0 &=  \ph{{} -\tht_{i0}\w\om_0 -\tht_{ij}}
      \hbox to 0pt{ \hss ${} - \tht_{0j}$}\w\om_j
         \ph{{} - I_{ijk}\,\tht_{k}\w\om_j}\,,\\
\d\om_i &=  {} -\tht_{i0}\w\om_0 -\tht_{ij}\w\om_j 
               - I_{ijk}\,\tht_{0k}\w\om_j\,.
\end{split}
\end{equation}

These equations have been derived starting with a generalized
Finsler structure~$(\Sigma,\iota)$.  It will be important in
what follows to know that there is the following sort of converse.
The proposition below may seem somewhat strange, but it is the
fundamental tool for ensuring that, once one has found
differential forms satisfying the appropriate structure equations,
they come from a (generalized) Finsler structure in a natural way.

\begin{proposition}\label{prop: recovery of gen Finler str}
Suppose that~$X$ is a manifold of dimension at least~$2n{+}1$
and that there exist linearly independent $1$-forms
$\om_0,\om_1,\dots,\om_n,\tht_{01},\dots,\tht_{0n}$
on~$X$ for which the equations~\eqref{eq: structure summary new}
hold for some~$1$-forms~$\theta_{ij}=-\theta_{ji}$ 
and some functions~$I_{ijk}=I_{jik}=I_{ikj}$ on~$X$.  

Suppose that there exist submersions~$\tau:X\to M^{n+1}$
and~$\psi:X\to\Sigma^{2n+1}$, respectively, with connected fibers 
such that the fibers of each are, respectively, 
the leaves of the Frobenius system 
generated by~$\{\om_0,\dots,\om_n\}$ 
and~$\{\om_0,\dots,\om_n,\tht_{01},\dots,\tht_{0n}\}$.

Then there exists an immersion~$\iota:\Sigma\to TM$ that defines
a generalized Finsler structure and a mapping~$\phi:X\to F$, where
$F$ is the canonical $\Or(n)$-structure associated
to~$(\Sigma,\iota)$, so that~$\phi$ pulls back the tautological
forms and connection forms on~$F$ to be the given forms on~$X$. 
\end{proposition}

\begin{proof} The key point is to explain how~$\iota$ is defined:
Take any vector field~$E$ on~$X$ (locally
defined, if necessary) for which~$\om_0(E)=1$ 
and~$\om_i(E)=\tht_{0i}(E)=0$.  Now define a mapping~$\tilde\iota:X\to TM$
by setting~$\tilde\iota(x) = \tau'\bigl(E(x)\bigr)$.  It is not 
difficult to show that~$\tilde\iota$ is constant on the fibers 
of~$\psi:X\to\Sigma$ and therefore that there exists a 
mapping~$\iota:\Sigma\to TM$ such that~$\tilde\iota = \iota\circ\psi$.
The mapping~$\phi$ is defined by a similar abstract diagram chase.

The remainder of the proof is a matter of checking details and can
be left to the reader.
\end{proof}

\subsubsection{A even wider sense of Finsler structure}
\label{sssec: wider sense}
At some point during this subsection (if not earlier),
the reader may have realized that even the generalization of
Finsler structure proposed in Definition~\ref{def: gen Finsler str}
is unnecessarily restrictive.  

The only ingredients used in the
construction are 
\begin{enumerate}
\item a $(2n{+}1)$-manifold~$\Sigma$, 
\item a contact form~$\om_0$, and 
\item a $\om_0$-Legendrian foliation~$\cM$ of~$\Sigma$ 
that satisfies the local convexity property needed to
ensure that the matrix~$H$ that shows up in~\eqref{eq: d om_i} is 
positive definite.  (This~$H$ represents a tensor
that is globally defined on~$\Sigma$ using only the data of~$\om_0$
and~$\cM$.)
\end{enumerate}

Thus, one could define a generalized Finsler structure to be 
a triple~$(\Sigma,\om_0,\cM)$ as above, subject to the appropriate
local convexity condition.  The construction of the canonical 
$\Or(n)$-structure ~$u:F\to\Sigma$ then proceeds just as before.

It will be useful to speak of generalized Finsler 
structures~$(\Sigma,\om_0,\cM)$
in this wider sense, so the reader should be alert for this usage in
the rest of this article.  The expression ``generalized Finsler structure
on~$M$'' will still be reserved for a pair~$(\Sigma,\iota)$ as
in~Definition~\ref{def: gen Finsler str}.

By Proposition~\ref{prop: recovery of gen Finler str}, 
any generalized Finsler structure in the wider sense is
locally realizable as a generalized Finsler structure on a manifold~$M$
and uniquely up to isometry to boot.  Thus, the extra generality is
only relevant when one does not have a manifold structure 
for the leaves for~$\cM$ explicitly in hand.

\subsection{The flag curvature}\label{ssec: flag curvature defined}
Differentiating the equations~\eqref{eq: structure summary new} 
(and reducing the second one modulo~$\om_1,\dots,\om_n$ to remove 
the derivatives of the functions~$I_{ijk}$) shows that there exist 
functions~$R_{0i0j}=R_{0j0i}$, $R_{0ijk}=-R_{0ikj}$, 
and $J_{ijk}=J_{jik}$ so that
\begin{equation}\label{eq: d theta 0i}
\d\tht_{0i} = -\tht_{0j}\w\tht_{ji} + R_{0i0j}\,\om_0\w\om_j
 + {\ts\frac12}R_{0ijk}\,\om_j\w\om_k + J_{ijk}\,\tht_{0k}\w\om_j\,.
\end{equation}
Moreover,~$R_{0ijk}+R_{0jki}+R_{0kij}=0$, just as in the Riemannian
case.  

Equation~\eqref{eq: d theta 0i} completes the first level 
of structure equations for the~$\Or(n)$-structure~$F$.

It will not be necessary to carry out a full development of the
structure equations here.  The interested reader is referred 
to~\cite{BCS} for a thorough treatment.

The most important aspect of these equations for the present
article is the so-called \emph{flag curvature}, represented
by the symmetric tensor
\begin{equation}\label{eq: flag curvature}
\cR = R_{0i0j}\,\,\om_i\circ\om_j\,.
\end{equation}

The geometric significance of the tensor~$\cR$ is that it furnishes
the lowest order term for the Jacobi equation that governs
the second variation of geodesics of the Finsler structure.  (This
should not be surprising since, as the structure equations show,~$\cR$ 
can be recovered from the \emph{second} Lie derivative of~$\gamma$ 
with respect to the Reeb vector field~$E$.)

\begin{definition}[Constant flag curvature]\label{def: cfc}
A generalized Finsler structure~$(\Sigma,\iota)$ is said to
have \emph{constant flag curvature} if there exists a constant~$c$
such that~$\cR = c\,\gamma$.
\end{definition}

\subsubsection{A PDE system}\label{sssec: a pde system}
While a local hypersurface~$\Sigma$ in~$TM$ depends essentially
on one arbitrary function of~$2n{+}1$ variables, 
the equation~$\cR=c\,\gamma$ is $\frac12\,n(n{+}1)$ 
non-linear partial differential equations for~$\Sigma$, thought
of as such a hypersurface. Thus, the condition of constant flag curvature 
is an overdetermined system of PDE for generalized Finsler structures
as soon as~$n>1$.

The generality of the solutions of this system up to 
local isometry is easily understood when~$n=1$ 
(see~\S\ref{ssec: dim 2} for a description in the case~$c=1$), 
but for~$n>1$ this system is overdetermined 
and it is not at all clear how many solutions there are, 
even locally.  In~\S\ref{sec: generality}, 
this question will be addressed for~$n\ge2$.

\subsubsection{The Jacobi equation}\label{sssec: jacobi}
When a generalized Finsler structure~$(\Sigma,\iota)$ has 
constant flag curvature~$c$, the Jacobi equation along any
of its geodesics is metrically conjugate to the Jacobi equation
of a geodesic in a Riemannian manifold of constant sectional
curvature~$c$.  Indeed, a Riemannian metric has constant
flag curvature when regarded as a Finsler structure if and
only if it has constant sectional curvature.

\subsubsection{Examples}\label{sssec: CFC examples}
Many (non-Riemannian) Finsler metrics
that have constant flag curvature are now known.  Here
are a few examples:
\begin{enumerate}
\item Hilbert showed that there is a canonically defined, complete 
Finsler metric of constant negative flag curvature on any bounded 
convex domain~$M$ in~$\bbR^{n+1}$ with strictly convex smooth boundary.
The geodesics of Hilbert's metric are the (open) line segments 
in~$M$.  Hence, Hilbert's examples are projectively flat.%
\footnote{
A Finsler metric~$\Sigma\subset TM$ is projectively flat 
if every point of~$M$ has a neighborhood~$U$ that can be embedded 
into~$\R{n+1}$ in such a way that it carries the $\Sigma$-geodesics in~$U$
to straight line segments.  In contrast to the Riemannian case,
where, by a theorem of Bonnet, projective flatness is equivalent to having 
constant sectional curvature, there are projectively flat 
Finsler metrics that do not have constant flag curvature and
Finsler metrics with constant flag curvature that are not 
projectively flat.}
Hilbert's Finsler metric is Riemannian if and only 
if the boundary of~$M$ is an ellipsoid.
\item  When~$n=2$, it is known~\cite{Br1} that the generalized Finsler 
metrics of constant flag curvature~$c$ depend essentially on two arbitrary
functions of two variables, up to local isometry.
Moreover, for~$c=1$, there are many examples defined globally on~$S^2$, 
at least as many (in some sense) as there are Zoll metrics on~$S^2$
(see~\S\ref{sssec: CFC 2-spheres}).
There is even a two-parameter family of mutually inequivalent, 
projectively flat Finsler metrics of constant flag curvature~$c=1$ 
on the $2$-sphere~\cite{Br1}.
\item Recently, Bao and Shen~\cite{BS} have constructed a one-parameter 
family of mutually inequivalent, left-invariant Finsler metrics on~$S^3$ 
that have constant flag curvature~$c=1$ that are not
projectively flat.
\end{enumerate}

However, the present knowledge of existence and/or classification
(in either the local or global case) of (generalized) Finsler structures 
with constant flag curvature is still very preliminary.  
The remaining sections of this article discuss various different aspects 
of this problem.

\section[A K\"ahler Structure]{A K\"ahler Structure}
\label{sec: kahler structure}

This section will be concerned with properties of the space of geodesics
of a generalized Finsler structure.

\begin{definition}[Geodesic simplicity]\label{def: geod simpl}
A generalized Finsler structure~$(\Sigma,\iota)$ on~$M^{n+1}$
is said to be \emph{geodesically simple} if the set~$Q$ of 
integral curves of its Reeb vector field~$E$ can be given the
structure of a smooth, Hausdorff manifold of dimension~$2n$
in such a way that the natural mapping~$\ell:\Sigma\to Q$ is
a smooth submersion.
\end{definition}

Any generalized Finsler structure is locally geodesically
simple, so for local calculations, one can always assume that~$\Sigma$
is geodesically simple.  Thus, throughout this section, unless it
is explicitly stated otherwise, it will be assumed that~$\Sigma$ is
geodesically simple.

Now, it is a classical fact drawn from the calculus of variations
that, for any generalized Finsler structure~$\Sigma$, the space~$Q$ inherits 
a canonical symplectic structure.  In fact, the $2$-form~$\d\om_0$ 
is manifestly invariant under the flow of~$E$ and satisfies~$E\lhk\d\om_0=0$, 
so there exists a symplectic form~$\Omega$ on~$Q$ so that~
\begin{equation}\label{eq: Omega defined}
\ell^*\Omega = -\d\om_0 = \tht_{0i}\w\om_i\,.
\end{equation}  
(The minus sign is for later convenience.)

In general, there is no natural metric on~$Q$.  However, when~$\Sigma$ 
has constant flag curvature, such a `metric'
does exist:

\begin{proposition}[Quotient quadratic form]\label{prop: quot quad form}
Suppose that~$(\Sigma,\iota)$ is a generalized Finsler structure with
constant flag curvature~$c$ that is geodesically simple.  Then there
exists a quadratic form~$\d\sigma^2$ on~$Q$ for which
\begin{equation}\label{eq: quad form}
\ell^*\bigl(d\sigma^2\bigr) = {\tht_{01}}^2+\dots+{\tht_{0n}}^2 
                               + c\,{\om_1}^2+\dots + c\,{\om_n}^2\,.
\end{equation}
\end{proposition}  

\begin{proof}  It suffices to show that the quadratic form on
the right hand side of~\eqref{eq: quad form} is invariant under
the flow of~$E$.  However, by hypothesis~$R_{0i0j} = c\,\delta_{ij}$.
Substituting this into the structure equations allows one to compute
the Lie derivative with respect to~$E$ of the right hand side 
of~\eqref{eq: quad form} and see that it is equal to zero.
\end{proof}

\subsection{The geodesic flow}\label{ssec: geodesic flow}
Assume that the generalized Finsler structure~$(\Sigma,\iota)$
has constant flag curvature~$c$.  Then the structure equations
derived so far take the form%
\footnote{ The reader may wonder about the choice of notation,
which is not entirely standard.  The letters~$I$ and~$J$ and
their signs were chosen to conform, 
insofar as possible, to Cartan's article~\cite{Ca1},
in which he treated the case~$n=1$.  This notation has turned
out to be very felicitous for computations in the constant
flag curvature case.}
\begin{equation}\label{eq: CFC c structure preliminary}
\begin{split}
\d\om_0 &=  \ph{\tht_{0j}\w\om_0 -\tht_{ij}}
             \hbox to 0pt{ \hss ${} - \tht_{0j}$}\w\om_j\,\\
\d\om_i &=  \tht_{0i}\w\om_0 -\tht_{ij}\w\om_j 
               - I_{ijk}\,\tht_{0k}\w\om_j\,,\\
\d\tht_{0i} &=  c\,\om_0\w\om_i - \tht_{ij}\w\tht_{0j} 
 + {\ts\frac12}R_{0ijk}\,\om_j\w\om_k + J_{ijk}\,\tht_{0k}\w\om_j\,,
\end{split}
\end{equation}
where the forms~$\tht_{ij}$ and functions~$R_{0ijk}$, $I_{ijk}$,
and~$J_{ijk}$ have the symmetries already discussed.

In what follows, the Lie derivative of a form or function with respect
to the Reeb vector field~$E$ will be denoted by an overdot.%
\footnote{ Bear in mind,
too, that~$E$ (which was originally defined on~$\Sigma$) has a canonical
lifting to the $\Or(n)$-bundle~$F$ so that it satisfies~$\om_0(E)=1$
while~$\om_i(E) = \tht_{0i}(E) = \tht_{ij}(E) = 0$.  No separate
notation will be introduced for the lifted vector field on~$F$.  The
reader should have no difficulty determining which is meant from
context.}

\begin{proposition}\label{prop: Reeb formulae}
The structure forms on the $\Or(n)$-bundle 
of a generalized Finsler structure of constant flag curvature~$c$ 
satisfy the following identities:
\begin{enumerate}
\item ${\dot\om}_i = -\tht_{0i}$ and ${\dot\tht}_{0i} = c\,\om_i$.
\item ${\dot I}_{ijk} = J_{ijk}$ and ${\dot J}_{ijk} = -c\,I_{ijk}$ 
\item $R_{0ijk} = 0$
\item ${\dot\tht}_{ij} = 0$.
\end{enumerate}
{\upshape(}In particular, $J$ is fully symmetric in its indices.{\upshape)}
\end{proposition}

\begin{proof}
The formulae in Item~$(1)$ follow immediately from the definition
of Lie derivative, the defining properties of~$E$, and the second
and third equations of~\eqref{eq: CFC c structure preliminary}. 

Now compute the Lie derivative with respect to~$E$ of the second line 
of~\eqref{eq: CFC c structure preliminary}, using the fact that this
operation is a derivation that commutes with exterior derivative, and
add the result to the third equation.  The result is
\begin{equation}\label{eq: E der one}
0 = -{\dot\tht}_{ij}\w\om_j + {\ts\frac12}R_{0ijk}\,\om_j\w\om_k
    + (J_{ijk}-{\dot I}_{ijk})\,\tht_{0k}\w\om_j\,.
\end{equation}
(The reader who performs this calculation will note that it uses
the fact that $I_{ijk}$ is fully symmetric in its indices.)  Since
$\tht_{ij}=-\tht_{ji}$ while~$I_{ijk}=I_{jik}$ and~$J_{ijk}=J_{jik}$,
it follows that~$J_{ijk}-{\dot I}_{ijk} = 0$ and thus that
\begin{equation}\label{eq: E der one simp}
0 = -{\dot\tht}_{ij}\w\om_j + {\ts\frac12}R_{0ijk}\,\om_j\w\om_k\,.
\end{equation}

In particular, the first equation of Item~$(2)$ is verified, which
shows that~$J$ is indeed fully symmetric in all its indices.  

Equation~\eqref{eq: E der one simp} also implies 
that there must exist (unique) functions~$T_{ijk}=-T_{jik}$ that satisfy
\begin{equation}\label{eq: E der one conseq}
{\dot\tht}_{ij} = -T_{ijk}\,\om_k \,,
\qquad\text{and}\qquad
R_{0ijk} = T_{ijk}-T_{ikj}\,.
\end{equation}

Now compute the Lie derivative with respect to~$E$ of the third line 
of~\eqref{eq: CFC c structure preliminary}, using the fact that this
operation is a derivation that commutes with exterior derivative, 
and subtract $c$ times the second equation from the result.  This 
yields the relation
\begin{equation}\label{eq: E der two}
0 = \bigl(2\,T_{ijk}-T_{ikj}-{\dot J}_{ijk}-c\,I_{ijk}\bigr)\,\om_k\w\tht_{0j} 
      + {\ts\frac12}{\dot R}_{0ijk}\,\om_j\w\om_k\,.
\end{equation}

Of course, this implies both~${\dot R}_{0ijk}=0$ and
\begin{equation}\label{eq: E der two simp}
0 = 2\,T_{ijk}-T_{ikj}-{\dot J}_{ijk}-c\,I_{ijk}\,.
\end{equation}
However, the symmetry of~$I$ and~$J$ and the 
skewsymmetry~$T_{ijk}=-T_{jik}$ now combine to show that
\begin{equation}\label{eq: E der two simp conseq}
T_{ijk} = 0
\qquad\qquad\text{and}\qquad\qquad
{\dot J}_{ijk}= -c\,I_{ijk}\,.
\end{equation}
This gives the second equation of Item~$(2)$ and, 
in view of~\eqref{eq: E der one conseq}, Items~$(3)$ and~$(4)$ as well.
\end{proof}

By Proposition~\ref{prop: Reeb formulae}, the structure 
equations for a generalized Finsler structure of constant
flag curvature~$c$ simplify to:

\begin{equation}\label{eq: CFC c structure final}
\begin{split}
\d\om_0 &=   - \tht_{0j}\w\om_j\,\\[5pt]
    \d\om_i &=  -\om_0\w\tht_{0i} -\tht_{ij}\w\om_j 
               - I_{ijk}\,\tht_{0k}\w\om_j\,,\\
\d\tht_{0i} &=   c\,\om_0\w\om_i - \tht_{ij}\w\tht_{0j} 
               + J_{ijk}\,\tht_{0k}\w\om_j\,,
\end{split}
\end{equation}
where $I$ and~$J$ are fully symmetric in their indices.
The similarity of the second and third lines is very 
suggestive and will be exploited in the next subsection.

\subsection{The K\"ahler structure}\label{ssec: Kahler structure}
The main concern of this article is the case of constant positive
flag curvature.  To treat this case, it suffices (by homothety)
to treat the case~$c=1$, so assume this from now on.

The pieces are now in place for the main result of this section:

\begin{theorem}\label{thm: Kahler str on Q}
Let~$(\Sigma,\iota)$ be a geodesically simple generalized Finsler
structure with constant flag curvature~$+1$.  Then the metric~$\d\sigma^2$
induced on~$Q$ is a K\"ahler metric, with K\"ahler form~$\Omega$.
Moreover,~$\Sigma$ has a natural immersion into the unit canonical
bundle of~$Q$ that commutes with the submersion~$\ell:\Sigma\to Q$.
\end{theorem}

\begin{proof} 
Define complex valued~$1$-forms
\begin{equation}\label{eq: zeta defined}
\zeta_i = \om_i - \iC\,\tht_{0i}\,.
\end{equation}
Using this notation and the condition~$c=1$, 
one finds that the pullbacks of~$\Omega$ and~$\d\sigma^2$
can be written in the form
\begin{equation}\label{eq: Omega and dsigma2 complex}
\begin{split}
\ell^*(\Omega) 
  &=  {\ts\frac\iC2}\,\zeta_i\w\ov{\zeta_i}\,,\\
\ell^*(\d\sigma^2) 
  &=  \zeta_1\circ\ov{\zeta_1} + \dots +  \zeta_n\circ\ov{\zeta_n}\,.
\end{split}
\end{equation}
Thus, the metric and $2$-form are algebraically compatible
and the $2$-form is closed.  The only condition remaining to verify
in order to show that this metric is K\"ahler with the given $2$-form 
as K\"ahler form is whether or not the almost complex structure
defined by this pair is integrable.

Now, the almost complex structure induced on~$Q$ is the one for which
the $\ell$-pullback of a $(1,0)$-form is a linear combination
of~$\zeta_1,\dots,\zeta_n$.  By the Newlander-Nirenberg Theorem, 
the integrability of the almost complex structure is equivalent 
to the condition that these forms define a differentially closed ideal.

To see this, note that the second and third 
structure equations~\eqref{eq: CFC c structure final} 
(when~$c=1$) can be combined in complex form as
\begin{equation}\label{eq: CFC 1 structure complex}
\d\zeta_i =  -\iC\,\om_0\w\zeta_i -\tht_{ij}\w\zeta_j 
  + {\ts\frac\iC2}(I_{ijk}+\iC\,J_{ijk})\,\ov{\zeta_j}\w\zeta_k\,.
\end{equation}
Thus, the complex Pfaffian system spanned by
the~$\zeta_i$ is Frobenius, as desired.

For the final statement, consider the complex-valued~$n$-form
defined on~$\Sigma$ by
\begin{equation}\label{eq: zeta vol def}
\zeta = \zeta_1\w\zeta_2\w\dots\w\zeta_n\,.
\end{equation}
Let~$K$ be the canonical bundle of~$Q$ (regarded as a complex manifold),
i.e., $K$ is the top exterior power of the complex cotangent bundle of~$Q$.
Let~$\Upsilon$ be the tautological holomorphic $n$-form on~$K$
and let $\Sigma(K)\subset K$ denote the circle bundle of unit complex
volume forms on~$Q$ with respect to the K\"ahler structure constructed
in the first part of the proof.  

Evidently, there is a unique smooth mapping~$\hat\ell:\Sigma\to\Sigma(K)$
that lifts~$\ell$ and satisfies~${\hat\ell}^*(\Upsilon) = \zeta$.  Because
$\zeta$ satisfies
\begin{equation}\label{eq: d zeta }
\d\zeta = -i\rho\w\zeta\,,\qquad\text{where}\qquad
\rho = n\,\om_0  -\Re\bigl((I_{iij}{-}\iC\,J_{iij})\,\zeta_j\bigr)\,,
\end{equation}
it follows that~$\hat\ell$ is an immersion.
\end{proof}

\subsection{A canonical flow}\label{ssec: canonical flow}
Proposition~\ref{prop: Reeb formulae} has another very interesting
consequence:  Let~$\Phi_t$ denote the (locally defined) flow of the 
Reeb vector field~$E$, lifted to~$F$.  Then 
Proposition~\ref{prop: Reeb formulae} implies (bearing in mind 
that~$c=1$ throughout this discussion):
\begin{equation}\label{eq: Phi_t on coframing}
\begin{split}
\Phi_t^*\om_0 &= \om_0\,,\\
\Phi_t^*\om_i &= \cos t\,\,\om_i - \sin t\,\,\tht_{0i}\,,\\
\Phi_t^*\tht_{0i} &= \sin t\,\,\om_i + \cos t\,\,\tht_{0i}\,,\\
\Phi_t^*\tht_{ij} &= \tht_{ij}\,.
\end{split}
\end{equation}
In particular, this flow commutes with the $\Or(n)$-action
that defines the bundle structures

Of course, $\Phi_t$ is only locally defined unless~$\Sigma$ 
is geodesically complete.  However, this does not prevent one
from defining a circle 
of coframings~$\bigl(\om^t_0,\om^t_i,\tht^t_{0i},\tht^t_{ij}\bigr)$ 
on the~$\Or(n)$-bundle~$F$ by the formulae
\begin{equation}\label{eq: circle of coframings}
\begin{split}
\om^t_0 &= \om_0\,,\\
\om^t_i &= \cos t\,\,\om_i - \sin t\,\,\tht_{0i}\,,\\
\tht^t_{0i} &= \sin t\,\,\om_i + \cos t\,\,\tht_{0i}\,,\\
\tht^t_{ij} &= \tht_{ij}\,.
\end{split}
\end{equation}
These forms (together with the appropriately rotated functions~$I_{ijk}$
and~$J_{ijk}$) evidently satisfy the structure 
equations~\eqref{eq: CFC c structure final}
with~$c=1$ for any value of~$t$ and so make~$F$ into the~$\Or(n)$-bundle
of a circle of generalized Finsler structures with 
constant flag curvature~$1$.  

The members of this circle of generalized Finsler structures are
generally not isometric among themselves.  Thus, this
constructs a nontrivial flow on the space of generalized Finsler
structures with constant flag curvature~$1$.  In the projectively
flat case (see the next section), however, the resulting generalized
Finsler structure is a fixed point of this flow.

\subsection{Related $G$-structures}\label{ssec: S1-On structure}
In the language of~$G$-structures, a K\"ahler structure
on a $2n$-manifold~$Q$ is the same thing 
as a torsion-free~$\Un(n)$-structure on~$Q$ where~$\Un(n)$
is embedded into~$\GL(2n,\bbR)$ in the usual way.  
Theorem~\ref{thm: Kahler str on Q} describes how a generalized
Finsler structure of constant flag curvature~$1$ determines
a natural K\"ahler structure on~$Q$.  

\subsubsection{An $S^1{\cdot}\Or(n)$-structure}\label{sssec: S1dotO(n) str}
Now, Proposition~\ref{prop: Reeb formulae} actually
implies that the structure on~$\Sigma$ determines
a canonical~$S^1{\cdot}\Or(n)$-structure on~$Q$, where
$S^1{\cdot}\Or(n)$ is the subgroup of~$\Un(n)$ generated
by~$\Or(n)$ and the central subgroup~$S^1\subset\Un(n)$
consisting of scalar multiplication by unit complex numbers.

This $S^1{\cdot}\Or(n)$-structure is defined as follows:
Let $\tau:F\to Q$ be the composition of the 
submersions~$u:F\to\Sigma$ and $\ell:\Sigma\to Q$.  
For~$f\in F$, 
the kernel of~$\tau'(f):T_fF\to T_{\tau(f)}Q$ is defined by
the equations~$\om_i=\tht_{0i}=0$.  It follows that
there is a unique isomorphism~$v(f):T_{\tau(f)}Q\to\C{n}$ 
for which~$\zeta = (\om_i-\iC\tht_{0i})$ satisfies
\begin{equation}\label{eq: v_f defined}
\zeta(w) = v(f)\bigl(\tau'(f)(w)\bigr)
\end{equation} 
for all~$w\in T_fF$.
Thus,~$v(f)$ is a $\C{n}$-valued coframe at~$\tau(f)$.  

Proposition~\ref{prop: Reeb formulae} 
implies that if~$\tau(f_1)=\tau(f_2)$, then $v(f_1)$
and~$v(f_2)$ are related by $v(f_1)=\e^{is}A\,v(f_2)$
for some~$s\in\bbR$ and~$A\in\Or(n)$. Thus, the fibers
of~$\tau$ are mapped into $S^1{\cdot}\Or(n)$-orbits 
in the bundle of $\C{n}$-valued coframes on~$Q$ and
dimension count plus the structure equations on~$F$
imply that~$v$ is actually a local diffeomorphism
on each fiber. Thus, there is a well-defined
$S^1{\cdot}\Or(n)$-structure~$\hat u:\hat F\to Q$ 
so that~$v:F\to \hat F$ immerses~$F$ as an open set
in~$\hat F$.

Equation~\eqref{eq: CFC 1 structure complex} can now
be interpreted as the structure equation on~$\hat F$,
where~$\om_0$ and the~$\theta_{ij}$ are the connection
forms for~$\hat F$ as a $S^1{\cdot}\Or(n)$-structure
over~$Q$. Note that, when~$n>1$, 
this $S^1{\cdot}\Or(n)$-structure has torsion 
unless~$I=J=0$. For the situation when~$n=1$, 
see~\S\ref{ssec: dim 2}.

\subsubsection{Extending~$\Sigma$}\label{sssec: extending Sigma}
Note that the quotient~$\hat\Sigma=\hat F/\Or(n)$
is naturally a circle bundle over~$Q$ and there
is a canonical map~$\Sigma\to\hat\Sigma$ that
is a local diffeomorphism, carrying the vector
field~$E$ into the infinitesimal generator of
the $S^1$-action.  Moreover, $\hat\Sigma$ carries
a natural generalized Finsler structure in the
wider sense of~\S\ref{sssec: wider sense}. 
The needed data are as follows:  First~$\om_0$
is well-defined on~$\hat\Sigma$ as the connection
form associated to the $S^1$-action.  Second,
there is a~$\cM$ of~$\hat\Sigma$ defined as the
image of the leaves of~$\om_0=\Re(\zeta)=0$ on~$\hat F$.  
One checks directly that this foliation extends
the foliation~$\cM$ on~$\Sigma$ and satisfies
the convexity assumptions as discussed 
in~\S\ref{sssec: wider sense}.  In particular,
it follows that, any geodesically simple
generalized Finsler structure~$\Sigma$ 
of constant flag curvature~$1$ can be 
canonically immersed in a generalized Finsler
structure~$\hat\Sigma$ that is geodesically complete,
with all geodesics closed of period~$2\pi$.

The leaves of~$\cM$ on~$\hat\Sigma$ project
via~$\ell$ to become Lagrangian submanifolds of~$Q$.
In fact, the (circle) fiber in~$\hat\Sigma$ over a 
point~$q\in Q$ represents a circle of Lagrangian
$n$-planes in~$T_qQ$ and the images of the leaves
of~$\cM$ are the Lagrangian submanifolds of~$Q$ 
whose tangent planes belong to~$\hat\Sigma$ when
regarded as a subset of the space of Lagrangian
planes of~$Q$.

Unfortunately, it can happen that the foliation~$\cM$
on~$\hat\Sigma$ is not simple.

\subsubsection{An $S^1{\cdot}\GL(n,\bbR)$-structure}
\label{sssec: S1dotGLnR str}
Finally, it is important to note that, although~$\hat F$
has torsion, it underlies a canonical
$S^1{\cdot}\GL(n,\bbR)$-structure on~$Q$ that does not:
Going back to~\eqref{eq: CFC 1 structure complex}, 
and setting
\begin{equation}\label{eq: sigma defined}
\sigma_{ij} = 
 {\ts\frac\iC2}(I_{ijk}-\iC\,J_{ijk})\,\zeta_k
-{\ts\frac\iC2}(I_{ijk}+\iC\,J_{ijk})\,\ov{\zeta_k}
=\ov{\sigma_{ij}}\,,
\end{equation}
one notes that, because~$I$ and~$J$ are symmetric
in all their indices,
\eqref{eq: CFC 1 structure complex}
can be written in the form
\begin{equation}\label{eq: s1 glnR str eqs}
d\zeta = -(\iC\,\om_0\,\I_n +\theta + \sigma)\w\zeta
\end{equation}
where~$\theta=(\theta_{ij})$ is real-valued and skewsymmetric 
while $\sigma=(\sigma_{ij})$ is real-valued and symmetric.
Since $(\iC\,\om_0\,\I_n +\theta + \sigma)$ takes values in
the Lie algebra of~$S^1{\cdot}\GL(n,\bbR)\subset \GL(2n,\bbR)$, 
it follows that, as promised, 
the $S^1{\cdot}\GL(n,\bbR)$-structure underlying~$\hat F$
has vanishing torsion.%
\footnote{ The reader will note that, when~$n>1$, 
a $S^1{\cdot}\GL(n,\bbR)$-structure has at most one
torsion-free connection.  Again, the case~$n=1$ is
special, but is treated by other means anyway.}

This is surprising because the 
subgroup $S^1{\cdot}\GL(n,\bbR)\subset\GL(2n,\bbR)$
acts irreducibly on~$\R{2n}$ and yet does not appear
on the accepted list~\cite{MS1,MS2} of irreducibly-acting holonomies
of torsion-free connections.  More will be said about
this point in~\S\ref{ssec: str eqs}, where
the structure equations will be investigated more thoroughly.

\section[Complex Hypersurfaces]
{Complex Hypersurfaces in $\bbC\bbP^{n+1}$}
\label{sec: cmplx hypersurf}

The goal of this section is to explain how
certain generalized Finsler structures on~$\bbR\bbP^{n+1}$ 
of constant flag curvature~$1$ can be constructed 
from complex hypersurfaces in~$\bbC\bbP^{n+1}$ 
that satisfy some genericity assumptions.

As usual, regard~$\bbR\bbP^{n+1} = \bbP^{n+1}$ as the space of real lines
in~$\R{n+2}$ through the origin.  The notation~$\tbbP^{n+1}$ will
denote its nontrivial double cover, i.e., the space of real \emph{rays}
in~$\R{n+2}$ emanating from the origin, or, equivalently, the set 
of oriented real lines in~$\R{n+2}$ through the origin.  If~$v\in\R{n+2}$
is nonzero, then~$[v]\in\bbP^{n+1}$ denotes the line spanned by~$v$
and~$[v]_+\in\tbbP^{n+1}$ denotes the ray spanned by~$v$.

Each oriented $2$-plane~$E\subset\R{n+2}$ through the origin 
determines an oriented line in~$\tbbP^{n+1}$ that is a closed circle, 
as follows:  If~$(v,w)$ is an oriented basis of~$E$, then the 
curve~$\gamma_{(v,w)}(s) = \bigl[(\cos s)\,v + (\sin s)\,w\bigr]_+$ 
is an oriented embedding of the circle into~$\tbbP^{n+1}$.  
It is easy to see that this oriented line (as an image) depends
only on the oriented $2$-plane~$E$ and not on the specific
oriented basis (more will be said about this below).  Thus,
$\Gr\!^+_2\bigl(\R{n+2}\bigr)$ parametrizes a family of oriented lines 
in~$\tbbP^{n+1}$ that has the property that there is a unique
such line passing through a given point and having a given
oriented tangent direction there.

A generalized Finsler structure~$(\Sigma,\iota)$ on~$\tbbP^{n+1}$ 
will be said to be \emph{rectilinear} if each of its (oriented) 
geodesics is a line (up to reparametrization).  (Note that
this is stronger than requiring saying that the
generalized Finsler structure be projectively flat.  It requires
that the geodesics actually \emph{be} lines in~$\tbbP^{n+1}$, 
not just that one can transform them into lines by local
reparametrizations in~$\tbbP^{n+1}$.)  In such a case, there
is a canonical submersion~$\lambda:\Sigma\to\Gr\!^+_2\bigl(\R{n+2}\bigr)$
whose fibers are unions of integral curves of~$E$.  (It need
not be true, in general, that the $\lambda$-fibers are connected.)

It was first observed by Funk~\cite{Fu1} that 
if~$D\subset\R{2}\subset\tbbP^2$ is a domain with 
a rectilinear Finsler structure~$\Sigma\subset TD$ 
whose flag curvature is identically~$1$, then every $\Sigma$-geodesic 
in~$D$ has a \emph{unit speed} parametrization of the form
\begin{equation}\label{eq: gamma unit speed}
\gamma_{(v,w)}(s) = \bigl[(\cos s)\,v + (\sin s)\,w\bigr]_+
\end{equation}
for some linearly independent~$v,w\in\R3$.  

Now, it is immediate that $\gamma_{(v_1,w_1)}$ and~$\gamma_{(v_2,w_2)}$
as defined by~\eqref{eq: gamma unit speed}
parametrize the same oriented geodesic with the same speed
if and only if~$[\![v_1+\iC\,w_1]\!] = [\![v_2+\iC\,w_2]\!]$
as points in~$\bbC\bbP^{n+1}\setminus\bbP^{n+1}$. 

Consequently, for a rectilinear Finsler structure~$\Sigma\subset TD$
with flag curvature identically~$1$, there is a refined `geodesic
map'~$\ell:\Sigma\to\bbC\bbP^2\setminus\bbP^2$ such that~$\Sigma$
is essentially recoverable from its image under~$\ell$.  (One
gets not only the geodesics, but their unit speed parametrizations
as well; hence one can recover their unit tangent vectors, i.e.,
$\Sigma\subset TD$.

Busemann~\cite{Bu} later pointed out that Funk's observation applies
just as well to a rectilinear Finsler structure of constant
flag curvature~$1$ on a projective space of any dimension.  In fact,
the argument is purely local, so that it applies to any rectilinear 
generalized Finsler structure~$(\Sigma,\iota)$ on~$\tbbP^{n+1}$  
with constant flag curvature~$1$.  The result is that one
has a natural mapping~$\ell:\Sigma\to\bbC\bbP^{n+1}\setminus\bbP^{n+1}$
for such a~$(\Sigma,\iota)$ whose differential has constant rank~$2n$
and whose fibers are discrete unions of integral curves of~$E$.  

Back in the case~$n=1$, Funk~\cite{Fu2} eventually observed that,
when~$\Sigma$ is a rectilinear Finsler structure on~$D\subset\tbbP^2$
with constant flag curvature~$1$, 
the $2$-dimensional image~$\ell(\Sigma)\subset\bbC\bbP^2\setminus\bbP^2$
is actually a holomorphic curve.  He then showed how, conversely, 
starting with a holomorphic curve in~$\bbC\bbP^2\setminus\bbP^2$ 
satisfying some open conditions, one could construct a (generalized)
rectilinear Finsler structure on a domain in~$\tbbP^2$ 
with constant flag curvature~$1$.  

In~\cite{Br3}, it was shown how Funk's construction could be 
globalized so as to classify the rectilinear Finsler structures 
on~$\tbbP^2$ with constant flag curvature~$1$.  It
was shown there that, up to projective equivalence, these structures
form a non-compact, $2$-parameter family.

The goal of this section is to explain how these constructions 
generalize to higher dimensions.  This turns out to be straightforward.  
However, it will be useful to examine the proofs directly via 
the moving frame for use in the next section.

\subsection{Some notation}\label{ssec: some notation}
It will be necessary to consider projective 
spaces of vector spaces with real or complex coefficients.
It will also sometimes be necessary to consider, 
for a real or complex vector space~$V$, 
the set~$\tbbP(V)$ of \emph{real rays} in~$V$, i.e., the equivalence
classes in~$V\setminus\{0\}$ generated by scalar multiplication
by positive real numbers.  For a nonzero vector~$v\in V$,
the notation~$[v]_+$ will be used for the real ray containing~$v$.
The standard notation~$[v]$ will be used for the real line containing~$v$
and, in case~$V$ is a complex vector space, the notation~$[\![v]\!]$
will be used for the complex line containing~$v$.

When~$V$ is a real vector space, $\tbbP(V)$ is naturally
the (non-trivial) double cover of~$\bbP(V)$ and is diffeomorphic 
to a sphere.  The space~$\tbbP\bigl(\bbR^{m}\bigr)$ will be 
denoted~$\tbbP^{m-1}$.

When~$V$ is a complex vector space, $\tbbP(V)$ is naturally
an $S^1$-bundle over~$\bbP(V)$ (which, as usual, denotes the
complex projectivization), and the natural mapping will be
denoted by~$\ell:\tbbP(V)\to\bbP(V)$.

The constructions will be designed so as to be equivariant
under the action of~$\SL(n{+}2,\bbR)$, so it will be useful
to consider some of the spaces on which this group acts.

First of all, there is~$\bbS = \tbbP(\C{n+2})\setminus\tbbP(\R{n+2})$.
A typical element of~$\bbS$ is of the form~$[v+\iC\,w]_+$ where~$v$
and~$w$ are linearly independent in~$\R{n+2}$.

There is a natural mapping~$\iota:\bbS\to T\/\tbbP^{n+1}$ that
sends~$[v+\iC\,w]_+\in\bbS$ to the velocity at~$t=0$ of the
curve~$\gamma_{(v,w)}:\bbR\to\tbbP^{n+1}$ defined
by~$\gamma_{(v,w)}(t) = [v + t\,w]_+$. 

The image of~$\ell:\bbS\to \bbC\bbP^{n+1}$ 
is the open set~$X_{n+1}=\bbC\bbP^{n+1}\setminus\bbR\bbP^{n+1}$.

The following diagrams may help to fix these homogeneous
spaces of~$\SL(n{+}2,\bbR)$ and $\SL(n{+}2,\bbR)$-equivariant 
mappings in their proper perspective:
\begin{equation}\label{eq: SL diagram}
\begin{CD}
\bbS @>{\iota}>> T\,\tbbP^{n+1}\\
@V{\ell}VV @VV{\pi}V\\
X_{n+1} @. \tbbP^{n+1}\\
@V{\lambda}VV @.\\
\Gr\!^+_2\bigl(\R{n+2}\bigr)\\
\end{CD}
\qquad\qquad
\begin{CD}
[v+\iC\,w]_+ @>{\iota}>> 
\left.\frac{d\hfil}{dt}\right|_{t=0}\left([v+t\,w]_+\right)\\
@V{\ell}VV @VV{\pi}V\\
[\![v+\iC\,w]\!] @. [v]_+\\
@V{\lambda}VV @.\\
[v\w w]_+\\
\end{CD}\hss
\end{equation}

The map~$\lambda$ is a surjective submersion 
and its fibers in~$X_{n+1}$ are Poincar\'e disks.  
In fact, the fiber over~$[v\w w]_+$ 
is one of the two disks cut out of the~$\bbC\bbP^1$ that 
is the projectivization of the (complex) span of~$v$ and~$w$ by
the removal of the real points.  The points of the fiber can be
thought of as metric structures on the oriented~$\bbR\bbP^1$ that
was removed.%
\footnote{ The reader is reminded that the foliation of~$X_{n+1}$ 
by the fibers of~$\lambda$ is \emph{not} holomorphic.  If it were,
then~$\Gr(2,n{+}2)$ would have a $\SL(n{+}2,\bbR)$-invariant
holomorphic structure, which it does not.}  

For further details about these and similar constructions,
see~\cite[\S4.1]{Br3}.

\subsection{The moving frame}\label{ssec: moving frame}
A moving frame for the homogeneous
space~$X_{n+1}$ can be constructed by starting with the standard
moving frame
\begin{equation}\label{eq: e frame}
\eb = (\eb_0\ \eb_1\ \cdots\ \eb_{n+1}):\SL(n{+}2,\bbR) 
   \longrightarrow  M_{n+2}(\bbR)
\end{equation}
with its canonical matrix-valued $1$-form~$\om=(\om^\beta_\alpha)$ 
satisfying the usual structure equations
\begin{equation}\label{eq: e str eqs}
\d\eb = \eb\,\om\,,\qquad\qquad
\d\om = { } -\om\w\om\,,\qquad\qquad
\text{and}\qquad\tr(\om)=0,
\end{equation}
and making a complex change of basis
\begin{equation}\label{eq: e to f}
(\fb_0\ \fb_1\ \cdots\ \fb_n) = (\eb_0\ \eb_1\ \cdots\ \eb_n)
\left[\begin{matrix} 
1 & 0 & 1\\
0 & \I_{n} & 0 \\
\iC & 0 & -\iC
\end{matrix}\right]\,.
\end{equation}
Note that~$\ov{\fb_0} = \fb_{n+1}$, $\ov{\fb_{n+1}} = \fb_0$, 
and $\ov{\fb_i} = \fb_i$ for~$1\le i\le n$.
The structure equation~$\d\fb = \fb \varphi$ 
can now be expanded in the following form (which establishes notation%
\footnote{ The introduction of the $\frac12$ coefficients simplifies
later normalizations.}
for the entries of~$\varphi$)
\begin{equation}\label{eq: f structure}
\d\left(\fb_0\quad\fb_i\quad\fb_{n+1}\right) 
= \left(\fb_0\quad\fb_j\quad\fb_{n+1}\right)
\left[\begin{matrix} 
\alpha   & -\frac12\,\ov{\pi_i}\  & \ov{\zeta^{n+1}}\ \\[5pt]
\zeta^j & \om^j_i  & \ov{\zeta^j}\\[5pt]
\zeta^{n+1} & -\frac12\,\pi_i  &  \ov{\alpha}
\end{matrix}\right].
\end{equation}
Note the trace relation
\begin{equation}\label{eq: f str trace relation}
\tr(\varphi) = \alpha + \ov{\alpha} + \om^i_i = 0
\end{equation}
and the second structure equation~$\d\varphi = -\varphi\w\varphi$,
which expands in the obvious way to provide formulae for~$\d\alpha$, etc.

By construction, the map~$[\fb_0]:\SL(n{+}2,\bbR)\to X_{n+1}$
is a surjective submersion and pulls back $(1,0)$-forms on~$X_{n+1}$
to be linear combinations of~$\{\zeta^1,\dots,\zeta^{n+1}\}$.

One more fact about this moving frame construction will be important.
The mapping
\begin{equation}\label{eq: framed incidence}
\left([\fb_0], [\fb_0\w\fb_1\w\dots\w\fb_n]\,\right):\SL(n{+}2,\bbR)
\longrightarrow X_{n+1}\times \left(\bbC\bbP^{n+1}\right)^*
\end{equation}
has, as its image, the set of pairs~$(p,H)$ where~$p\in X_{n+1}$
is any point and~$H\subset\bbC\bbP^{n+1}$ is any complex hyperplane 
through~$p$ that is transverse to the $\lambda$-fiber through~$p$.
The easy verification of this fact is left to the reader.

\subsection{Transverse, convex hypersurfaces}\label{ssec: convex hyps}
Now let~$Q\subset\bbC\bbP^{n+1}$ be a (not necessarily compact) 
nonsingular complex hypersurface that is transverse to the fibers of~$\lambda$.

Let~$\Sigma_Q$ denote the preimage~$\ell^{-1}(Q)\subset\bbS$, and
let~$\iota_Q:\Sigma_Q\to T\/\tbbP^{n+1}$ be the restriction of~$\iota$
to~$\Sigma_Q$.  Then~$\Sigma_Q$ is a smooth manifold of dimension~$2n{+}1$
and it is not hard to see that the assumption that~$Q$ be transverse 
to the fibers of~$\lambda$ implies that 
the map~$\iota_Q:\Sigma_Q\to T\/\tbbP^{n+1}$
is a radially transverse immersion.

Define the \emph{first order frame bundle}~$F^1_Q\subset\SL(n{+}2,\bbR)$
of~$Q$ to be the set of~$\fb$ in~$\SL(n{+}2,\bbR)$ so that~$[\fb_0]$ 
lies in~$Q$ and the hyperplane~$[\fb_0\w\fb_1\w\dots\w\fb_n]$ 
is tangent to~$Q$ at~$\fb_0$.  The projection~$[\fb_0]:F^1_Q\to Q$ 
is then a bundle over~$Q$ whose fibers are left cosets of the 
subgroup~$G_1\subset\SL(n{+}2,\bbR)$ consisting of the
set of matrices of the form 
\begin{equation}\label{eq: G1 element}
\left[\begin{matrix} 
1 & 0 & 1\\
0 & \I_{n} & 0 \\
\iC & 0 & -\iC
\end{matrix}\right]\,
\left[\begin{matrix} 
a & 0 & 0\\
0 & A & 0 \\
0 & 0 & \ov{a}
\end{matrix}\right]\,
\left[\begin{matrix} 
1 & 0 & 1\\
0 & \I_{n} & 0 \\
\iC & 0 & -\iC
\end{matrix}\right]^{-1}\,.
\end{equation}
where~$A\in\GL(n,\bbR)$ and~$a\in\C*$ satisfy~$a\ov{a}\det(A)=1$.

Pulling the forms on~$\SL(n{+}2,\bbR)$ back to~$F^1_Q$,
one has the relations
\begin{equation}
\zeta^{n+1} = 0,\qquad\qquad\qquad
\zeta^1\w\cdots\w\zeta^n\not=0.
\end{equation}
Since the structure equations entail
\begin{equation}
\d\zeta^{n+1} 
= (\alpha-\ov{\alpha})\w\zeta^{n+1}+{\ts\frac12}\,\pi_j\w\zeta^j\,,
\end{equation}
it follows that~$\pi_j\w\zeta^j = 0$.  By Cartan's Lemma, 
there exist functions~$H_{ij} = H_{ji}$ on~$F^1_Q$ so that
\begin{equation}
\pi_i = H_{ij}\,\zeta^j\,.
\end{equation}
The structure equations entail
\begin{equation}\label{eq: dalpha}
\d\alpha = {\ts\frac12}\,\ov{\pi_i}\w\zeta^i 
         = -{\ts\frac12}\,\ov{H_{ij}}\,\zeta^i\w\ov{\zeta^j}
\end{equation}
and
\begin{equation}\label{eq: dzeta a}
\d\zeta^i = (\delta^i_j\,\alpha - \om^i_j)\w\zeta^j\,.
\end{equation}
After an application of Cartan's Lemma, these relations yield
\begin{equation}\label{eq: dH}
\d H_{ij} - (\alpha+\ov{\alpha})\,H_{ij} 
          + H_{ik}\,\om^k_j + H_{kj}\,\om^k_i   = H_{ijk}\,\zeta^k
\end{equation}
for some functions~$H_{ijk} = H_{jik} = H_{ikj}$ on~$F^1_Q$.

In particular, the complex-valued quadratic 
form~$\cH = H_{ij}\,\zeta^i\,\ov{\zeta^j}$
is well-defined on~$Q$, as are its real and imaginary parts.

\begin{example}[Standard null quadric]\label{ex: standard null quad}
When restricted to~$\SO(n{+}2)\subset\SL(n{+}2,\bbR)$,
the map~$[\fb_0]$ has image equal to the standard null quadric, 
whose homogeneous equation is
\begin{equation}\label{eq: null quad}
(z^0)^2 + (z^1)^2 + \cdots + (z^{n+1})^2 = 0.
\end{equation}
In this case, the structure matrix~$\varphi$ reduces to
\begin{equation}\label{eq: varphi on null quad}
\left[\begin{matrix} 
\iC\,\om^0_{n+1}&-\frac12\,\ov{\zeta_i}\ &\ov{\zeta^{n+1}}\ \\[5pt]
\zeta^j &\om^j_i & \ov{\zeta^j}\\[5pt]
\zeta^{n+1} & -\frac12\,\zeta_i  &  -\iC\,\om^0_{n+1}
\end{matrix}\right],
\qquad\qquad\text{where}\qquad\qquad
\om^i_j + \om^j_i = 0.
\end{equation}
In particular,~$H=\I_{n}$ in this case.
\end{example}

With this example in mind, the following definition will be adopted:

\begin{definition}\label{def: convex}
A smooth embedded complex hypersurface~$Q\subset X_{n+1}$ that is
transverse to the~$\lambda$-fibers will be said to be \emph{convex}
if the function $\Re(H)$ on~$F^1_Q$ 
takes values in positive definite matrices, or, equivalently,
if the quadratic form~$\Re(\cH)$ is positive definite on~$Q$.
\end{definition}

It is an elementary exercise to check that the condition that~$Q$
be convex is equivalent to the condition that~$(\Sigma_Q,\iota_Q)$
satisfies the local convexity condition given 
in Definition~\ref{def: gen Finsler str} that is needed to ensure 
that it be a generalized Finsler structure on~$\tbbP^{n+1}$.
Unfortunately, it need not satisfy the fiber-connectedness 
hypotheses given in~Definition~\ref{def: gen Finsler str},
so this is only a generalized Finsler structure in the wider
sense.

The way is now prepared for stating the main result of this
section, which generalizes the construction in~\cite{Br3} for~$n=1$
that was based on an idea of Funk~\cite{Fu2}.

\begin{theorem}\label{thm: convex gives CFC}
Suppose that the complex hypersurface~$Q\subset X_{n+1}$ 
is transverse and convex. Then the generalized Finsler 
structure~$(\Sigma_Q,\iota_Q)$ on~$\tbbP^{n+1}$
is rectilinear and has constant flag curvature~$+1$.

Conversely, for any rectilinear generalized Finsler structure~$(\Sigma,\iota)$
on~$\tbbP^{n+1}$ with constant flag curvature~$+1$, 
the canonical geodesic map~$\ell:\Sigma\to X_{n+1}$
has, as its image, a complex (immersed) hypersurface~$Q$
that is transverse and convex.
\end{theorem}

The proof of the first half of the theorem will be given 
in the following subsection.  (Afterwards, the proof of the
converse statement can safely be left to the reader.)
Note that the example of the standard null quadric, 
Example~\ref{ex: standard null quad},
must correspond to the Riemannian metric of constant curvature~$1$
on~$\tbbP^{n+1}\simeq S^{n+1}$.  
A further discussion of examples will be taken up after the proof.

\subsection{Structure reduction}\label{ssec: structure reduc}
Assume for the rest of this section that~$Q$ is convex.  

The equation~$\Re(H)=\I_{n}$ defines a sub-bundle~$F^2_Q\subset F^1_Q$
whose structure group is the group~$G_2\subset G_1$ consisting
of the matrices of the form~\eqref{eq: G1 element} 
with~$a\ov{a}=1$ and~$A\in\SO(n)$.  Thus $G_2\simeq S^1\times\SO(n)$.  
Henceforth, all functions and forms will be regarded as pulled 
back to~$F^2_Q$, though, as is customary in moving frame calculations, 
this pullback will not be notated.

It will be useful to separate~$\zeta^i$ into its real and imaginary
parts, so introduce  real-valued forms~$\om_i$ and~$\tht_{0i}$ 
for~$1\le i\le n$ by the equations%
\footnote{ This notation is chosen so as to agree with the 
notation in earlier sections.}
\begin{equation}\label{eq: zeta re and im}
\zeta^i = \om_i -\iC\,\tht_{0i}\,.
\end{equation}

Now the equation
\begin{equation}\label{eq: H reduced}
H = \I_{n}+\iC\,Y
\end{equation}
holds, where~$Y$ is symmetric and real-valued.  
Define $1$-forms~$\rho$ and~$\om_0$ 
so as to separate~$\alpha$ into its real and imaginary parts as 
\begin{equation}\label{eq: alpha reduced}
\alpha = \rho - \iC\,\om_0\,.
\end{equation}
Then~\eqref{eq: dalpha} becomes
\begin{equation}\label{eq: dalpha reduced}
\d\bigl(\rho-\iC\,\om_0\bigr) 
= -{\ts\frac12}\,(\delta_{ij}-\iC\,Y_{ij})\,\zeta^i\w\ov{\zeta^j}.
\end{equation}
Separating this equation into its real and imaginary parts yields
\begin{equation}\label{eq: d rho}
\d\rho = {\ts\frac\iC2}\,Y_{ij}\,\zeta^i\w\ov{\zeta^j}
       = Y_{ij}\,\tht_{0i}\w\om_j
\end{equation}
and
\begin{equation}\label{eq: d omega0}
\d\om_0 = -{\ts\frac\iC2}\,\zeta^i\w\ov{\zeta^i} = -\tht_{0i}\w\om_i\,.
\end{equation}
Now \eqref{eq: dH} can be written in the form
\begin{equation}\label{eq: dH reduced}
\d H_{ij} =  (\delta_{ik}\,\rho-\om^k_i)\,H_{kj}
           + (\delta_{kj}\,\rho-\om^k_j)\,H_{ik}
           + H_{ijk}\,\zeta^k.
\end{equation}
It will be useful to separate this into its real and imaginary parts.
First, set
\begin{equation}\label{eq: H_ijk re and im}
H_{ijk} = 2\,(J_{ijk}+\iC\,I_{ijk}),
\end{equation}
where~$I_{ijk}$ and~$J_{ijk}$ are real-valued
and then define new $1$-forms~$\tht_{ij}=-\tht_{ji}$ 
and~$\sigma_{ij}=\sigma_{ji}$ by the relations
\begin{equation}\label{eq: w_ij and sigma_ij defined}
\delta_{ij}\,\rho-\om^i_j = {} - \tht_{ij} - \sigma_{ij}\,.
\end{equation}
The real part of~\eqref{eq: dH reduced} can now be written in the form
\begin{equation}\label{eq: sigma_ij}
\sigma_{ij} = {\ts\frac12}\Re\bigl(H_{ijk}\,\zeta^k\bigr)
            = J_{ijk}\,\om_k + I_{ijk}\,\tht_{0k}\,.
\end{equation}

Moreover, the structure 
equation~$\d\zeta = \bigl(\delta^i_j\alpha - \om^i_j\bigr)\w\zeta^j$
separates into real and imaginary parts as
\begin{equation}\label{eq: om a and om 0a}
\begin{split}
\d\om_{i\ph0} 
  &=\tht_{0i}\w\om_0-\tht_{ij}\w\om_{j\ph0}-I_{ijk}\,\tht_{0k}\w\om_j\,,\\
\d\tht_{0i} 
  &=\om_{0}\w\om_{i\ph{0}}-\tht_{ij}\w\tht_{0j} +J_{ijk}\,\tht_{0k}\w\om_j\,,
\end{split}
\end{equation}

The reader will recognize equations~\eqref{eq: d omega0} 
and~\eqref{eq: om a and om 0a} as the structure equations of 
the canonical~$\SO(n)$-bundle of a generalized Finsler structure
of constant flag curvature~$1$.  

Of course, there needs to be a base manifold of dimension~$n{+}1$,
but this is easily constructed:  Note that, by the structure equations
and definitions so far
\begin{equation}\label{eq: de_0}
\d\eb_0 = \d\bigl(\Re(\fb_0)\bigr) 
= \eb_0\,\rho + \eb_i\,\om_i + \eb_{n+1}\,\om_0\,.
\end{equation}
Thus, the mapping~$[\eb_0]_+:F^2_Q\to\tbbP^{n+1}$ is a submersion
and its fibers are (unions of) leaves of the 
Frobenius system~$\om_0=\om_1=\cdots=\om_n=0$.  
Moreover, the fibers of the map~$[\fb_0]_+:F^2_Q\to\bbS$
are the~$\SO(n)$-orbits on~$F^2_Q$ and the image of this
map is~$\Sigma_Q$, by definition.  

It now follows from the structure equations 
that~$\bigl(\Sigma_Q,\iota_Q\bigr)$ is
a rectilinear generalized Finsler structure on~$\tbbP^{n+1}$ 
with constant flag curvature~$+1$, as desired.

\subsection{Examples}\label{ssec: Q examples}
It is now time to consider some examples of transverse, convex
hypersurfaces. 

\begin{example}[Non-real Hyperquadrics]\label{ex: non-real quadrics}
Let~$Q\subset X_{n+1}$ be a hypersurface so that
the induced generalized Finsler structure is actually
a Finsler structure on~$\tbbP^{n+1}$.  By construction, this 
means that~$Q$ is compact and hence algebraic.  Moreover, since
each geodesic in~$\tbbP^{n+1}$ occurs with two orientations, it
follows that~$Q$ must meet each $\lambda$-fiber transversely in
two points.  It follows that~$Q\subset\bbC\bbP^{n+1}$ has degree two,
i.e., is a hyperquadric and has no real points.

Now, a hyperquadric~$Q$ with no real points is $\SL(n{+}2,\bbR)$-equivalent
to a unique hyperquadric of the form
\begin{equation}\label{eq: p parameter quadric}
(z^0)^2 + e^{\iC p_1}\,(z^1)^2 +\cdots + e^{\iC p_{n+1}}\,(z^{n+1})^2 =0,
\end{equation}
where~$p = (p_1,\,\dots,\,p_{n+1})$ is a real vector satisfying
\begin{equation}\label{eq: p inequalities}
0=p_0\le p_1\le\cdots\le p_{n+1}<\pi\,.
\end{equation}

Conversely, it is not difficult to see that the quadric~$Q_p$
defined by~\eqref{eq: p parameter quadric} where the~$p_i$
are subject to~\eqref{eq: p inequalities} is both transverse
and convex.  Moreover, it is easy to see that distinct values of~$p$
give rise to non-isometric Finsler structures.

Thus, this provides an $(n{+}1)$-parameter family of distinct,
rectilinear Finsler structures with constant
flag curvature~$1$ on~$S^{n+1} = \tbbP^{n+1}$.

Only the case~$p = (0,\dots,0)$ is Riemannian.  When the $p_i$ 
(including~$p_0=0$)
are distinct, the group of isometries of the corresponding Finsler
metric is discrete, but it has positive dimension when two or more 
of the~$p_i$ are equal.
\end{example}

\begin{example}[Prescribed Indicatrix]\label{ex: prescribed indicatrix}
Theorem~\ref{thm: convex gives CFC} can be used to construct rectilinear
Finsler structures with constant flag curvature~$1$ 
and a prescribed tangent indicatrix at one point.  
In fact, one has the following result:

\begin{proposition}\label{prop: prescribed indicatrix}
Let~$[v]_+\in\tbbP^{n+1}$ be any point and let~$S\subset T_{[v]_+}\tbbP^{n+1}$
be a compact, real-analytic hypersurface that is strictly convex towards
the origin in~$T_{[v]_+}\tbbP^{n+1}$.  Then there is an open neighborhood~$U$
of~$[v]_+$ in~$\tbbP^{n+1}$, together with a Finsler structure~$(\Sigma,\iota)$
on~$U$, so that~$\iota(\Sigma)\subset TU$ contains~$S$ 
and so that~$(\Sigma,\iota)$
is rectilinear and has constant flag curvature~$1$.
\end{proposition}

\begin{proof}
Choose a hyperplane~$W\subset\bbR^{n+2}$ 
that is transverse to the line~$[v]$ 
and set~$\hat S 
= \left\{ [v+\iC\,w]_+\in\bbS
            \ \vrule\ w\in W,\ \iota\bigl([v+\iC\,w]_+\bigr)\in S\right\}$.
Of course,~$\hat S$ is diffeomorphic to~$S\simeq S^n$.  
The image~$\ell\bigl(\hat S\bigr)\subset X_{n+1}$ 
is a totally real, real analytic $n$-dimensional submanifold of~$X_{n+1}$ 
whose complexified tangent space is transverse to the fibers of~$\lambda$.
Thus, there exists a unique complex hypersurface~$Q\subset X_{n+1}$ that
contains~$\ell\bigl(\hat S\bigr)$.  By restricting~$Q$ to a sufficiently
small tubular neighborhood of~$\ell\bigl(\hat S\bigr)$ 
(in some metric on~$X_{n+1}$), one can assume that~$Q$ is embedded and
everywhere transverse to the fibers of~$\lambda$ (since it is
along~$\ell\bigl(\hat S\bigr)$.  Moreover, the hypothesis that~$S$ is
strictly convex towards the origin in~$T_{[v]_+}\tbbP^{n+1}$ implies
that~$Q$ is convex (in the sense of Definition~\ref{def: convex}) 
on a neighborhood of~$\ell\bigl(\hat S\bigr)$, so by shrinking~$Q$ again
if necessary, one can assume that~$Q$ is convex everywhere.  

Consider the corresponding~$\bigl(\Sigma_Q,\iota_Q\bigr)$, which
is a rectilinear generalized Finsler structure on~$\tbbP^{n+1}$ 
with constant flag curvature~$1$.  By construction, the
fiber~$\Sigma_{[v]_+}=\hat S$ is compact and convex.  It is now not 
difficult to see that there must be an open neighborhood~$U$ of~$[v]_+$
in~$\tbbP^{n+1}$ with the property that, for all~$[v']_+\in U$, 
the fiber~$\Sigma_{[v']_+}$ is also compact and nonempty.  This~$U$
is the desired neighborhood.
\end{proof}

\begin{remark}[A more general construction]\label{rem: gen tan indic}
Note that the argument in the proof does not construct all 
of the possible rectilinear Finsler structures on a neighborhood of~$[v]_+$ 
with constant flag curvature~$1$ and 
with the given tangent indicatrix at the point~$[v]_+$.  

In fact, if $\lambda:\hat S\to \bbR$ is any real analytic function,
set
\begin{equation}\label{eq: gen tan indic}
\hat S_\lambda = \left\{\ [v+\iC\,(w+\lambda\,v)]_+\in\bbS
            \ \vrule\ [v+\iC\,w]_+\in \hat S\ \right\}.
\end{equation}
Then one can use~$\ell\bigl(\hat S_\lambda\bigr)$ 
instead of~$\ell\bigl(\hat S\bigr)$ to generate a complex
hypersurface and the construction proceeds as before.  This more
general construction does give all of the the possible rectilinear 
Finsler structures 
on a neighborhood of~$[v]_+$ with constant flag curvature~$1$ 
with the given tangent indicatrix at the point~$[v]_+$.

Of course, these methods do not give any easy method to estimate
how large the domain~$U$ will be.

In some sense, this construction is the positive curvature
analog of Hilbert's construction of rectilinear Finsler metrics
with constant flag curvature~$-1$ on convex domains in~$\R{n+1}$.
\end{remark}

\end{example}

\section[Generality]{Generality}\label{sec: generality}

\subsection{The case of dimension~$2$}\label{ssec: dim 2}
For comparison, the local description of generalized Finsler
metrics on surfaces with constant flag curvature~$1$
will be recalled from~\cite{Br1}.

The structure equations in case~$n=1$ take the form
\begin{equation}\label{eq: om0 om1 tht1}
\begin{split}
\d\om_{0\ph0} &= -\tht_{01}\w\om_1\,,\\
\d\om_{1\ph0} 
  &=-\om_0\w\tht_{01} - I\,\tht_{01}\w\om_1
   = -(\om_0 - I\,\om_1 + J\,\tht_{01})\w\tht_{01} \,,\\
\d\tht_{01} 
  &=\ph{-}\om_{0}\w\om_{1\ph{0}} + J\,\tht_{01}\w\om_1
   = \ph{-}(\om_0 - I\,\om_1 + J\,\tht_{01})\w\om_1\,,
\end{split}
\end{equation}
where, throughout this subsection, $I_{111}$ and ~$J_{111}$ 
will be written as~$I$ and~$J$, respectively.
These are the structure equations
on the $\Or(1)$-structure~$F$ over~$\Sigma$.  By passing
to a double cover if necessary, it will be assumed that
these equations hold on~$\Sigma$ itself.  

Assuming that~$\Sigma$ is geodesically simple
with geodesic projection~$\ell:\Sigma\to Q$, 
Proposition~\ref{prop: Reeb formulae} implies that, 
not only do there exist a metric~$\d\sigma^2$
and area form~$\Omega$ on~$Q$ 
satisfying
\begin{equation}\label{eq: area form and metric}
\ell^*(\Omega) = \tht_{01}\w\om_1\,,
\qquad\qquad\qquad
\ell^*(\d\sigma^2) = {\tht_{01}}^2+{\om_1}^2\,,
\end{equation}
but there also exists a $1$-form~$\beta$ on~$Q$ satisfying
\begin{equation}\label{eq: extra 1-form}
\ell^*(\beta) = { } - I\,\om_1 + J\,\tht_{01}\,.
\end{equation}

A glance at~\eqref{eq: om0 om1 tht1} coupled with
knowledge of the structure equations of a Riemannian metric 
shows that
\begin{equation}\label{eq: d varphi}
\d\beta = (1-K)\,\Omega
\end{equation}
where~$K$ is the Gauss curvature of the metric~$\d\sigma^2$.

Conversely, suppose that one has a surface~$Q$ 
endowed with a metric~$\d\sigma^2$ with Gauss curvature~$K$, 
an area form~$\Omega$, and a $1$-form~$\beta$
that satisfies~$\d\beta = (1-K)\,\Omega$.  Let~$\ell:\Sigma\to Q$
be the oriented orthonormal frame bundle of~$Q$ endowed with
the metric~$\d\sigma^2$ and orientation~$\Omega$.  Then the
usual tautological and connection forms~$\eta_1,\eta_2,\eta_{12}$
defined on~$\Sigma$ satisfy
\begin{equation}\label{eq: Q pullbacks}
\ell^*(\d\sigma^2) = {\eta_1}^2+{\eta_2}^2\,,\qquad\quad
\ell^*(\Omega) = \eta_1\w\eta_2\,,\qquad\quad
\ell^*(\beta) = { } - I\,\eta_2 + J\,\eta_1\,,
\end{equation}
for some functions~$I$ and~$J$ on~$\Sigma$, the structure equations
\begin{equation}\label{eq: Q str eqs}
\d\eta_1 = -\eta_{12}\w\eta_2\,,\qquad\quad
\d\eta_2 =  \eta_{12}\w\eta_1\,,\qquad\quad
\d\eta_{12} = \ell^*(K)\,\eta_1\w\eta_2\,,
\end{equation}
and the equation
\begin{equation}\label{eq: d varphi on frame}
\d({}-I\,\eta_2 + J\,\eta_1) = \bigl( 1-\ell^*(K)\bigr)\,\eta_1\w\eta_2\,.
\end{equation}
Consequently, setting
\begin{equation}\label{eq: eta to Finsler str}
\om_0 = -\eta_{12} + I\,\eta_2 - J\,\eta_1\,,\qquad\quad
\om_1 =  \eta_2\,,\qquad\quad
\tht_{01} = \eta_1\,,
\end{equation}
yields a coframing on~$\Sigma$ that satisfies the
structure equations for a generalized Finsler structure
with constant flag curvature~$1$.

Thus, the local prescription for generalized Finsler surfaces
with constant flag curvature~$1$ is equivalent
to prescribing data on a surface: a metric~$d\sigma^2$,
its area form~$\Omega$, and a $1$-form~$\beta$ that satisfies
the equation~$\d\beta = (1-K)\,\Omega$.  Up to local
isometry, a metric~$\d\sigma^2$ on a surface depends 
on one arbitrary function of two variables and the 
$1$-form~$\beta$ is determined up to the addition of
an exact $1$-form~$\d f$, which is also one function of two
variables.  

Thus, (local) generalized Finsler structures 
for surfaces with constant flag curvature~$1$
depend on two arbitrary functions of two variables.

\subsubsection{$\beta$-geodesics}\label{sssec: beta geodesics}
Generally, given a metric~$\d\sigma^2$ with area form~$\Omega$ 
on a surface~$Q$ and a $1$-form~$\beta$, a curve~$\gamma\subset Q$
that satisfies~$\kappa_\gamma\,\d s = \beta_{|\gamma}$ 
will be called a \emph{$\beta$-geodesic} with respect 
to~$\bigl(\d\sigma^2,\Omega\bigr)$.  Here, $\kappa_\gamma$
represents the geodesic curvature of~$\gamma$ when one 
fixes an orientation of~$\gamma$.  Of course, reversing
the orientation of~$\gamma$ reverses the sign of both its
arc length~$\d s$ and its geodesic curvature~$\kappa_\gamma$,
so the expression~$\kappa_\gamma\,\d s$ is unchanged.

The orientation of the surface is significant:  
A curve~$\gamma$ is a $\beta$-geodesic 
with respect to~$\bigl(\d\sigma^2,\Omega\bigr)$ 
if and only if it is a~$(-\beta)$-geodesic 
with respect to~$\bigl(\d\sigma^2,-\Omega\bigr)$.

Just as in the case of ordinary geodesics (i.e., the $0$-geodesics), 
there is a unique $\beta$-geodesic with respect 
to~$\bigl(\d\sigma^2,\Omega\bigr)$
with any given initial point and direction on the surface~$Q$.

The $1$-form~$\beta$ is sometimes called the ``magnetic field''
for particles moving on~$Q$.

\subsubsection{CFC $2$-spheres}\label{sssec: CFC 2-spheres}
Now return to the case of a geodesically simple 
generalized Finsler structure~$\ell:\Sigma\to Q$
endowed with a coframing~$(\om_0,\om_1,\tht_{01})$ 
satisfying~\eqref{eq: om0 om1 tht1}.
Define~$\d\sigma^2$, $\Omega$, and~$\beta$ on~$Q$
by~\eqref{eq: Q pullbacks}.

The leaves of the system~$\om_0=\om_1=0$ on~$\Sigma$, i.e., 
the fibers of a realization~$\pi{\circ}\iota:\Sigma\to M$
as a generalized Finsler structure on a surface~$M^2$, 
are then mapped to the $\beta$-geodesics with respect 
to~$\bigl(\d\sigma^2,\Omega\bigr)$

For example, when~$\beta=0$, these curves are geodesics.
Of course, the condition~$\beta=0$ implies that~$K=1$,
so that these are just the geodesics on a standard $2$-sphere~$Q$ 
of constant Gauss curvature~$1$.  The corresponding~$M$ is just
the $2$-sphere of oriented geodesics on the standard $2$-sphere.
More interesting examples will be constructed below.

In general, if the data~$\bigl(\d\sigma^2,\Omega,\beta\bigr)$
on~$Q$ has the property that the~$\beta$-geodesics
with respect to~$\bigl(\d\sigma^2,\Omega\bigr)$ are all closed, 
then they lift to closed curves in~$\Sigma$ regarded as the 
unit tangent bundle of~$Q$ and the quotient surface~$M$ 
will exist globally.

There are now two elementary results to note.  Each is
a calculation that can be left to the reader.
\begin{proposition}\label{prop: conformal varphi geodesics}
Let~$Q$ be a surface endowed with a metric~$\d\sigma^2$,
an area form~$\Omega$, and a $1$-form~$\beta$. For any
function~$L>0$ on~$Q$ define
\begin{equation}\label{eq: conformal formulae}
\d\bar\sigma^2 = L\,\d\sigma^2\,,\qquad\quad
\bar\Omega = L\,\Omega\,,\qquad\quad
\bar\beta = \beta + *\d\bigl(\log\sqrt{L}\bigr)\,.
\end{equation}
Then the $\bar\beta$-geodesics with respect
to~$\bigl(\d\bar\sigma^2,\bar\Omega\bigr)$ 
are the same as the $\beta$-geodesics with respect
to~$\bigl(\d\sigma^2,\Omega\bigr)$. \qed
\end{proposition}
\begin{proposition}\label{prop: pos curv constr}
Let~$Q$ be a surface endowed with a metric~$\d\sigma^2$ with
Gauss curvature~$K>0$ and area form~$\Omega$.  Then the data
\begin{equation}\label{eq: pos curv constr}
\d\bar\sigma^2 = K\,\d\sigma^2\,,\qquad\quad
\bar\Omega = K\,\Omega\,,\qquad\quad
\bar\beta =  *\d\bigl(\log\sqrt{K}\bigr)\,,
\end{equation}
satisfy~$\d\bar\beta = (1-\bar K)\,\bar\Omega$, where~$\bar K$
is the Gauss curvature of~$\d\bar\sigma^2$. \qed
\end{proposition}

Recall that a metric~$\d\sigma^2$ on the $2$-sphere is said 
to be a \emph{Zoll metric} (see~\cite[Chapter 4]{Be}) 
if all of its geodesics are closed and
of length $2\pi$.  It is elementary to show that, in this case,
the space of oriented $\d\sigma^2$-geodesics 
is itself a $2$-sphere~$M$.

In resolving a question of Funk, Guillemin~\cite{Gu} has shown
that there exist many Zoll metrics near the metric of constant
Gauss curvature~$1$ on~$S^2$.  See~\cite[Chapter 4]{Be}, for
another account and further discussion of related problems.

\begin{theorem}\label{thm: Zoll construction}
Let~$\d\sigma_0^2$ be a Zoll metric on~$Q = S^2$ with positive Gauss
curvature~$K_0$.  Let~$\Omega_0$ be the area form of~$\d\sigma_0^2$ and
let~$M\simeq S^2$ be the space of oriented $\d\sigma_0^2$-geodesics on~$Q$.
Then there exists a unique Finsler structure~$\Sigma\subset TM$ 
on~$M$ with constant flag curvature~$1$ whose geodesic 
projection~$\ell:\Sigma\to Q$ has the induced data
\begin{equation}\label{eq: Zoll constr}
\d\sigma^2 = K_0\,\d\sigma_0^2\,,\qquad\quad
\Omega = K_0\,\Omega\,,\qquad\quad
\beta = *\d\bigl(\log\sqrt{K_0}\bigr)\,.
\end{equation}
\end{theorem} 

\begin{proof}
By hypothesis, the $0$-geodesics of~$\bigl(\d\sigma_0^2,\Omega_0\bigr)$
are all closed, so, by Proposition~\ref{prop: conformal varphi geodesics},
the $\beta$-geodesics of~$\bigl(\d\sigma^2,\Omega\bigr)$ (which are
the same) are also closed.  Moreover, 
by Proposition~\ref{prop: pos curv constr}, the 
data~$\bigl(\d\sigma^2,\Omega,\beta\bigr)$ 
satisfy~$\d\beta = (1-K)\,\Omega$ where~$K$ is the Gauss curvature 
of~$\d\sigma^2$.  By the discussion at the beginning of this subsection,
there is a canonically constructed coframing~$(\om_0,\om_1,\tht_{01})$ 
on~$\ell:\Sigma\to Q$, the unit tangent bundle of~$\d\sigma^2$ over~$Q$,
that satisfies the structure equations~\eqref{eq: om0 om1 tht1} of
a generalized Finsler structure of constant flag curvature~$1$ and
that induces the given data~$\bigl(\d\sigma^2,\Omega,\beta\bigr)$
on~$Q$, its space of geodesics.  Because its foliation~$\cM$ given
by~$\om_0=\om_1=0$ has closed leaves and, in fact, has~$M$ as its
leaf space, Proposition~\ref{prop: recovery of gen Finler str} shows
that there is an immersion~$\iota:\Sigma\to TM$ that realizes~$\Sigma$
as a generalized Finsler structure on~$M$.  The reader can easily check
that~$\Sigma$ is, in fact, an embedding and defines a genuine Finsler
structure on~$M$, as desired.
\end{proof}

\begin{remark}[Other global possibilities]\label{rem: other glob poss}
Theorem~\ref{thm: Zoll construction} provides one way to construct 
data~$\bigl(\d\sigma^2,\Omega,\beta\bigr)$ on~$S^2$ 
satisfying~$\d\beta = (1-K)\,\Omega$ and the condition 
that the~$\beta$-geodesics 
with respect to~$\bigl(\d\sigma^2,\Omega\bigr)$ be closed. 

Note that this Zoll construction 
only produces data~$\bigl(\d\sigma^2,\Omega,\beta\bigr)$
with~$\d(*\beta)=0$.  In fact, by writing~$*\beta = \d u$
for some function~$u$ (uniquely determined up to an additive constant), 
one can recover the original Zoll metric
from this data by dividing~$\d\sigma^2$ by~$e^{2u}$.
Thus, the Finsler structure~$\Sigma\subset TM$ determines
the original Zoll metric.  
Consequently, Theorem~\ref{thm: Zoll construction} provides 
an injection of the set
of isometry classes of Zoll metrics with positive Gauss curvature 
into the set of isometry classes of Finsler metrics on~$S^2$ with
constant flag curvature~$1$.  

The Zoll method is far from the only method of constructing global
examples, though it is the most general found so far.  
For example, one can find other examples by assuming rotational symmetry
in the data.  Also, the projectively flat examples
constructed in Example~\ref{ex: non-real quadrics} (with $n=1$) 
do not arise from the Zoll construction (except for the Riemannian one).

None of these examples (other than the Riemannian one) are reversible,
i.e., $\Sigma\not=-\Sigma\subset TM$.  In fact, the 
data~$\bigl(\d\sigma^2,\Omega,\beta\bigr)$ on~$Q$ give rise to a 
reversible Finsler structure on~$M$ if and only if there exists
a fixed-point free involution~$\iota:Q\to Q$ that fixes~$\d\sigma^2$
and reverses~$\Omega$ and~$\beta$.  No such example with~$\beta\not=0$
is known at present (nor has it been ruled out).
\end{remark}

\subsection{The structure equations in higher dimensions}
\label{ssec: str eqs}
As was already mentioned in~\S\ref{sssec: S1dotGLnR str}, 
a generalized Finsler structure~$\bigl(\Sigma,\iota\bigr)$ 
with constant flag curvature~$1$ that is geodesically simple 
induces a torsion-free $S^1{\cdot}\GL(n,\bbR)$-structure 
on the space~$Q$ of geodesics.  It turns out that this 
construction is essentially reversible, as will now be
explained.  Then, in later subsections, this reversibility 
will be used to investigate the generality of generalized
Finsler structures with constant flag curvature~$1$.

For the rest of this section, the assumption~$n>1$ will be
in force.  

\subsubsection{A circle of totally real $n$-planes}
\label{sssec: circle of n-planes}
Since~$S^1{\cdot}\GL(n,\bbR)$ is a subgroup of~$\GL(n,\bbC)$
(assuming their standard embeddings into~$\GL(2n,\bbR)$, 
a torsion-free $S^1{\cdot}\GL(n,\bbR)$ on a $2n$-manifold~$Q$
underlies an integrable almost complex structure.  Geometrically,
the reduction from an integrable almost complex structure
to an $S^1{\cdot}\GL(n,\bbR)$-structure is represented by
the choice of a totally real $n$-plane in each tangent
space, defined up to multiplication by~$\e^{\iC\theta}$.
Equivalently, one has a subbundle~$R\subset\Gr(n,TQ)$
of totally real tangent $n$-planes~$E\subset T_qQ$
(i.e., $E\cap \iC\,E = \{0_q\}$) for which the fiber
over each point~$R_q\subset R$ consists of the complex
multiples of single totally real $n$-plane.

Conversely, the choice of such a circle bundle~$R\subset\Gr(n,TQ)$
over a complex $n$-manifold~$Q$ defines a $S^1{\cdot}\GL(n,\bbR)$-structure
$q:F\to Q$: A coframing~$u:T_qQ\to\C{n}$ belongs to the structure~$F$ 
if and only if~$u$ carries the elements of the fiber~$R_q$ 
to the $n$-planes~$\e^{\iC\theta}\R{n}$.

Given such a circle bundle~$R\subset\Gr(n,TQ)$,
a~$n$-dimensional submanifold~$P\subset Q$ will be said
to belong to~$R$ if its tangent plane at every point is
an element of~$R$.  Belonging to~$R$ is an overdetermined system 
of first order partial differential equations for 
submanifolds~$P\subset Q$.  If~$P\subset Q$ belongs to~$R$, then
it has a canonical lifting~$\tau:P\to R$ defined by~$\tau(q)=T_qP$
for~$q\in P$.  This will be called the \emph{tangential lifting}
of~$P$.

It is easy to see that, 
for every~$n$-plane~$E\in R$, there is at most one connected
$n$-dimensional submanifold~$P\subset Q$ that belongs to~$R$
and has~$E$ as its tangent plane.  (This uses the hypothesis~$n>1$.)
The bundle~$R$ and, by association, 
the corresponding $S^1{\cdot}\GL(n,\bbR)$-structure~$q:F\to Q$
will be said to be \emph{integrable} if every element of~$R$
is tangent to an $n$-manifold that belongs to~$R$.  The condition
of being integrable is equivalent to the condition that~$R$ be
foliated by the tangential lifts of the $n$-manifolds that belong
to~$R$. 

\begin{example}[Generalized Finsler structures]\label{ex: gen Fin circle}
If~$(\Sigma,\iota)$ is a generalized Finsler structure on~$M^{n+1}$
with constant flag curvature~$1$ that is geodesically simple, with
geodesic projection~$q:\Sigma\to Q$, then the images~$q(\Sigma_x)\subset Q$
for~$x\in M$ belong to the canonical 
torsion-free~$S^1{\cdot}\GL(n,\bbR)$-structure
constructed in~\S\ref{sssec: S1dotGLnR str}.  Their liftings foliate an
open set in the associated circle bundle~$R$ and, in fact, $R$ is
integrable, as will be seen below.
\end{example}

\subsubsection{Torsion-free structures}
\label{sssec: torsion-free}
An $S^1{\cdot}\GL(n,\bbR)$-structure $q:F\to Q$ will be said
to be \emph{torsion-free} if it admits a connection without
torsion.  

Denote the Lie algebra of~$S^1{\cdot}\GL(n,\bbR)\subset \GL(2n,\bbR)$
by~$\eut\oplus\eugl(n,\bbR)\subset\eugl(2n,\bbR)$.  It is straightforward
to compute that the first prolongation%
\footnote{ See~\cite{Br2}
for information related to prolongation.} 
%~$\bigl(\eut\oplus\eugl(n,\bbR)\bigr)^{(1)}$
of this subalgebra of~$\eugl(2n,\bbR)$ vanishes (this uses
the assumption~$n>1$).  Consequently, if~$q:F\to Q$ does
admit a torsion-free connection, it admits only one. 

It will be necessary to examine the structure equations 
of~$F$ in the torsion-free case, in particular, to compute
the space of curvature tensors 
of torsion-free $S^1{\cdot}\GL(n,\bbR)$-structures.

Let~$\zeta = \left(\zeta^i\right)$ be the tautological $\C{n}$-valued $1$-form 
on~$F$. The assumption that~$F$ be torsion-free is equivalent to assuming
that there exist on~$F$ a $1$-form~$\om_0$ and a $\eugl(n,\bbR)$-valued
$1$-form~$\phi = \left(\phi^i_j\right)$ so that the 
\emph{first structure equation}
\begin{equation}\label{eq: F 1st str eq}
\d\zeta = -\bigl(\iC\,\om_0\,\I_n + \phi \bigr) \w\zeta
\end{equation}
holds.  These forms~$\om_0$ and~$\phi$ are the connection
forms of the structure.  

The \emph{second structure equation} will give
expressions for the curvature forms
\begin{equation}\label{eq: F curv expr}
\Omega_0 = \d\om_0\,,\qquad\qquad\qquad
\Phi = \d\phi + \phi \w\phi\,,
\end{equation}
that are based on the \emph{first Bianchi identity}
\begin{equation}\label{eq: F 1st Bianchi}
0 = -\bigl(\iC\,\Omega_0\,\I_n + \Phi \bigr) \w\zeta\,,
\end{equation}
which is derived by computing the exterior derivative 
of~\eqref{eq: F 1st str eq}.  This computation,
which is left to the reader, has the following result.

\begin{proposition}[Second structure equations]\label{prop: 1st Bianchi}
If~$n>2$, there exist on~$F$ real-valued functions~$b_{ij}=b_{ji}$
and~$r^i_{jkl}=r^i_{kjl}=r^i_{jlk}$ so that
\begin{equation}\label{eq: 1st Bianchi n>2}
\begin{split}
\Omega_0 & = -\iC\,b_{kl}\,\zeta^k\w\ov{\zeta^l}\,,\\
\Phi^i_j &= b_{jl}\,\bigl(\zeta^i\w\ov{\zeta^l}+\ov{\zeta^i}\w\zeta^l\,\bigr)
             + \iC\,r^i_{jkl}\,\zeta^k\w\ov{\zeta^l}\,.
\end{split}
\end{equation}
If~$n=2$, there exist on~$F$, in addition to the real-valued 
functions~$b_{ij}=b_{ji}$ and~$r^i_{jkl}=r^i_{kjl}=r^i_{jlk}$,
a complex-valued function~$A$ and a real-valued function~$a$
so that 
\begin{equation}\label{eq: 1st Bianchi n=2}
\begin{split}
\Omega_0 & = \Im\left(A\,\zeta^1\w\zeta^2\right)
              + 3a\,\bigl(\zeta^1\w\ov{\zeta^2}-\zeta^2\w\ov{\zeta^1}\,\bigr)
              - \iC\,b_{kl}\,\zeta^k\w\ov{\zeta^l}\,,\\
\Phi^i_j &= \delta^i_j\Re\left(A\,\zeta^1\w\zeta^2\right)
       + b_{jl}\,\bigl(\zeta^i\w\ov{\zeta^l}+\ov{\zeta^i}\w\zeta^l\,\bigr)
     + \iC\,(r^i_{jkl}+a(\epsilon_{jk}\delta^i_l{+}\epsilon_{jl}\delta^i_k))
                \,\zeta^k\w\ov{\zeta^l}\,.
\end{split}
\end{equation}
where~$\epsilon_{ij}=-\epsilon_{ji}$ and~$\epsilon_{12}=1$. \qed
\end{proposition}

\begin{remark}[Prolongation algebra]\label{rem: prolongation algebra}
Let~$V$ be an abstract real vector space of dimension~$n$ 
with complexification~$V^\bbC$.  The algebra~$\eugl(V)$ is naturally
included into~$\eugl(V^\bbC)$ and one can consider 
the Lie algebra~$\eug = \bbC{\cdot}\I_{V^\bbC} + \eugl(V)$ 
as a (real) sub-algebra 
of~$\eugl(V^\bbC)$. This is a proper subalgebra as long as~$n>1$.

It has already been remarked that, when~$n>1$, 
the first prolongation vanishes: $\eug^{(1)}=0$.  
Proposition~\ref{prop: 1st Bianchi} computes~$\cK(\eug)$, 
the space of curvature tensors of a torsion-free $\eug$-connection.  
The result is 
\begin{equation}\label{eq: K(g)}
\cK(\eug) = 
\begin{cases}
S^2(V^*)\oplus V{\ot}S^3(V^*),  & \text{when~$n>2$},\\
S^2(V^*)\oplus V{\ot}S^3(V^*)\oplus \bbC\oplus\bbR,  & \text{when~$n=2$}.
\end{cases}
\end{equation}
Note that the generic element in~$S^2(V^*)\oplus V{\ot}S^3(V^*)\subset 
\cK(\eug)$ does not lie in~$\cK(\euh)$ 
for any proper sub-algebra~$\euh\subset\eug$, so a $\eug$-connection
whose curvature assumes such a generic value will have holonomy equal
to the full group~$S^1{\cdot}\GL(n,\bbR)$.  Thus, Berger's first
criterion for~$S^1{\cdot}\GL(n,\bbR)$ to exist as the holonomy of
a torsion-free connection is satisfied.
\end{remark}

\begin{corollary}[Integrability]\label{cor: integrability}
When~$n>2$, a torsion-free~$S^1{\cdot}\GL(n,\bbR)$-structure
$q:F\to Q$ is integrable.  When~$n=2$, such a structure is integrable
if and only if the functions~$A$ and~$a$ vanish identically on~$F$.
\end{corollary}

\begin{proof}
The integrability condition is equivalent to the condition
that the Pfaffian system on~$F$ generated by~$\om_0$ 
and the components of~$\Im(\zeta)$ be Frobenius.  
By Proposition~\ref{prop: 1st Bianchi}, this condition is satisfied when~$n>3$
and is satisfied when~$n=2$ if and only if~$A=a=0$.
\end{proof}

Since the only $S^1{\cdot}\GL(n,\bbR)$-structures that
arise in the study of generalized Finsler structures with
constant flag curvature~$1$ are integrable and torsion-free, 
only the integrable, torsion-free case will be considered 
further in this article.  In order to have a uniform
notation, let~$\cK_\circ(\eug)\subset\cK(\eug)$ denote
the subspace consisting of the tensors of integrable, torsion-free 
$S^1{\cdot}\GL(n,\bbR)$-structures.  Thus~$\cK_\circ(\eug)\simeq
S^2(V^*)\oplus V{\ot}S^3(V^*)$ for all~$n\ge2$.

For an integrable, torsion-free $S^1{\cdot}\GL(n,\bbR)$-structure~$q:F\to Q$, 
the structure equations 
derived so far can be written in the form
\begin{equation}\label{eq: F 1st and 2nd str eq}
\begin{split}
\d\zeta^i &= -\bigl(\iC\,\delta^i_j\,\om_0 + \phi^i_j \bigr) \w\zeta^j\\
\d\om_0  & = -\iC\,b_{kl}\,\zeta^k\w\ov{\zeta^l}\,,\\
\d\phi^i_j+\phi^i_k\w\phi^k_j 
       &= b_{jl}\,\bigl(\zeta^i\w\ov{\zeta^l}+\ov{\zeta^i}\w\zeta^l\,\bigr)
             + \iC\,r^i_{jkl}\,\zeta^k\w\ov{\zeta^l}\,.
\end{split}
\end{equation}
where~$b_{ij}=b_{ji}$ and~$r^i_{jkl}=r^i_{kjl}=r^i_{jlk}$ are
real-valued functions on~$F$.

For later purposes, it will be necessary to understand 
the second Bianchi identity as well.  This is computed by applying the
exterior derivative to the second and third equations 
of~\eqref{eq: F 1st and 2nd str eq} and working out the consequences.  

	The result of the computation is that 
there exist \emph{complex}-valued functions~$B_{ijk}=B_{jik}=B_{ikj}$
and~$R^i_{jklm}=R^i_{kjlm}=R^i_{jlkm}=R^i_{jkml}$ on~$F$ so that
\begin{equation}\label{eq: F 2nd Bianchi}
\begin{split}
\d b_{ij} &= b_{kj}\phi^k_i+b_{ik}\phi^k_j 
+ \Re\left(B_{ijk}\zeta^k\right)\,,\\
\d r^i_{jkl} & = -r^m_{jkl}\phi^i_m
                  +r^i_{mkl}\phi^m_j+r^i_{jml}\phi^m_k+r^i_{jkm}\phi^m_l\\
       &\qquad\qquad
 +\Re\left(\left(R^i_{jklm}
         -\iC\,(\delta^i_j\,B_{klm}+\delta^i_k\,B_{ljm}
                      +\delta^i_l\,B_{kjm})\right)\zeta^m\right)\,.
\end{split}
\end{equation}

\begin{remark}[Prolongation algebra continued]
\label{rem: prolongation algebra cont}
In the notation of Remark~\ref{rem: prolongation algebra},
this calculation has the following interpretation:  
When~$n>2$, this second Bianchi identity calculation 
determines the space~$\cK^1(\eug)$, i.e., the space of 
first covariant derivatives of curvature tensors
of torsion-free~$S^1{\cdot}\GL(n,\bbR)$-structures.  
Then formula~\eqref{eq: F 2nd Bianchi}
implies the isomorphism
\begin{equation}\label{eq: K1(g)}
\cK^1(\eug) = S^3(V^*)^\bbC \oplus \left(V{\ot}S^4(V^*)\right)^\bbC.
\end{equation}
When~$n=2$, this is not the calculation of~$\cK^1(\eug)$ since
the integrability condition~$A=a=0$ has been imposed.  However,
in this case, the formula above does describe the space of
covariant derivatives of curvature tensors of \emph{integrable}
torsion-free~$S^1{\cdot}\GL(n,\bbR)$-structures, which, it turns
out, is the space that needed to be computed for applications in
this article anyway, since this space is the prolongation
of~$\cK_\circ(\eug)$ (regarded as a second-level tableau)
in either case.

In particular, it follows that~$\cK^1(\eug)\not=0$ for all~$n\ge 2$, 
so Berger's second criterion for $S^1{\cdot}\GL(n,\bbR)$ to be
the holonomy of a torsion-free connection that is not locally
symmetric is also satisfied.
\end{remark}

\begin{proposition}[Involutivity]\label{prop: invol tab}
Regard~$\cK_\circ(\eug)$ as a subspace 
of~$\eug \ot \Lambda^2\bigl((V^\bbC)^*\bigr)$, i.e., as
a second-level tableau. This subspace is involutive, 
with Cartan characters
\begin{equation}\label{eq: Cartan characters}
s_k = 
\begin{cases}
0, & k = 0,1,\\
k-1+n\bigl(n+(k{-}2)(n{+}1{-}k)\bigr), & 2\le k\le n{+}1,\\
0, & n{+}1 < k \le 2n.
\end{cases}
\end{equation}
The characteristic variety of this tableau consists of the
covectors~$\xi\in\bbP\bigl((\bbC\ot V)^*\bigr)\simeq\bbP^{2n-1}$
of the form~$\lambda\ot \xi'$ for~$\xi'\in\bbP(V^*)$ 
and is of degree~$n{+}1$ in~$\bbP\bigl((\bbC\ot V)^*\bigr)$.
\end{proposition}

\begin{proof}
This is a straightforward calculation:  Choose a flag that is
non-characteristic with respect to the claimed characteristic
variety.  One then finds that the characters of this flag are 
as given in~\eqref{eq: Cartan characters}.  
However, by combinatorics, one sees that, 
not only does one have the identity
\begin{equation}\label{eq: character sum}
s_2 + \dots + s_{n+1} 
= \binom{n+1}{2} + n\,\binom{n+2}{3}
= \dim \cK_\circ(\eug),
\end{equation}
but also that Cartan's test is satisfied, i.e., 
\begin{equation}\label{eq: Cartan's test}
2\,s_2 + \dots + (n{+}1)\,s_{n+1} 
= 2\,\binom{n+2}{3} + 2n\,\binom{n+3}{4}
= \dim \left(\cK_\circ(\eug)\right)^{(1)},
\end{equation}
as was verified in the computation of the second Bianchi identity
for integrable, torsion-free $S^1{\cdot}\GL(n,\bbR)$-structures.
\end{proof}

This has an immediate consequence:  

\begin{theorem}\label{thm: structure generality}
Up to diffeomorphism, the local, integrable,
torsion-free $S^1{\cdot}\GL(n,\bbR)$-structures
in dimension~$2n$ depend on~$n(n{+}1)$ functions
of~$n{+}1$ variables.  Moreover, for any 
curvature tensor in~$\cK_\circ(\eug)$, 
there exists an integrable, torsion-free 
$S^1{\cdot}\GL(n,\bbR)$-structure on a neighborhood
of~$0\in\R{2n}$ that assumes this value at~$0$.
\end{theorem}

\begin{proof}
These results follow from the usual Cartan-style construction
of an exterior differential system whose integral
manifolds are the the local, integrable,
torsion-free $S^1{\cdot}\GL(n,\bbR)$-structure plus
the algebraic result of Proposition~\ref{prop: invol tab}.
The proof is similar in all details to those 
executed in~\cite{Br2}, to which the reader is 
referred if more detail is needed.
\end{proof}

\begin{remark}[Exotic Holonomies]\label{rem: exotic holonomy}
In~\cite{Br2}, the existence of torsion-free 
$S^1{\cdot}\GL(2,\bbR)$-structures on $4$-manifolds whose
canonical connections have holonomy equal to $S^1{\cdot}\GL(2,\bbR)$
was established and in \cite{Br4}, the existence of torsion-free 
$H_\lambda{\cdot}\SL(2,\bbR)$-structures on $4$-manifolds whose
canonical connections have holonomy~$H_\lambda{\cdot}\SL(2,\bbR)$
for any $1$-parameter subgroup~$H_\lambda\subset\bbC^*$ (other than~$\bbR^*$)
was established. 

Using Proposition~\ref{prop: invol tab}, one can similarly demonstrate
the existence of torsion-free $H_\lambda{\cdot}\SL(n,\bbR)$-structures 
on $2n$-manifolds whose canonical connections have holonomy equal 
to~$H_\lambda{\cdot}\SL(n,\bbR)$ for any $1$-parameter 
subgroup~$H_\lambda\subset\bbC^*$ (other than~$\bbR^*$).  

This is interesting because these holonomy groups are not on Berger's
original lists of holonomies of irreducible holonomy torsion-free
connections (and hence fall into the category of `exotic' holonomies) 
and also were apparently overlooked in the recent classification 
of such holonomies by Merkulov and Schwachh\"ofer~\cite{MS1,MS2}.
\end{remark}

All that remains is to tie the geometry of these
structures to that of generalized Finsler structures
with constant flag curvature~$1$.  The key to doing
this is structure reduction.

Note that, if~$q:F\to Q$ is a torsion-free~$S^1{\cdot}\GL(n,\bbR)$-structure,
then the curvature $2$-form~$\Omega_0$ is actually the $q$-pullback
of a $2$-form that is well-defined on~$Q$.  By abuse of notation,
the symbol~$\Omega_0$ will be used to denote this $2$-form on~$Q$
as well.  If the structure is also integrable, then, by 
Proposition~\ref{prop: 1st Bianchi} and Corollary~\ref{cor: integrability}, 
the form $\Omega_0$ is of type~$(1,1)$ on~$Q$.  

Say that the structure~$F$ is \emph{positive} 
if the symmetric matrix~$b = (b_{ij})$
takes values in positive definite matrices or, 
equivalently, if~$-\Omega_0$ is a positive $(1,1)$-form on~$Q$, i.e.,
it defines a K\"ahler structure on~$Q$. 
In this case, there is a canonical~$S^1{\cdot}\Or(n)$-substructure
$F_\circ\subset F$ defined by the equations~$b_{ij} = \frac12\delta_{ij}$.
This will be called the \emph{K\"ahler reduction} of the 
torsion-free, integrable $S^1{\cdot}\GL(n,\bbR)$-structure~$F$.

Now the preparations have been made for the statement of the
final result of this article:

\begin{theorem}\label{thm: CFC1 realization}
Let~$q:F\to Q$ be a torsion-free 
$S^1{\cdot}\GL(n,\bbR)$-structure, assumed integrable 
if~$n=2$. If~$-\Omega_0$ is a positive~$(1,1)$-form
on~$Q$, then the K\"ahler reduction of~$F$ defines a 
generalized Finsler structure with constant flag
curvature~$1$.
\end{theorem}

\begin{proof}
This is a matter of computation and expansion of the definitions.
The point is that if one reduces to the locus in~$F$ where
$b_{ij} = \frac12\delta_{ij}$, this clearly defines an 
$S^1{\cdot}\Or(n)$-substructure $F_\circ\subset F$ as mentioned
above.  One can then write~$\zeta_i = \om_i-\iC\,\tht_{0i}$
and write~$\phi^i_j = \theta_{ij} +\sigma_{ij}$,  just
as in the previous section.  Then the first equation 
of~\eqref{eq: F 2nd Bianchi} shows how one can define~$I_{ijk}$
and~$J_{ijk}$ in terms of the real and imaginary parts of~$B_{ijk}$
so that equations~\eqref{eq: CFC c structure final} hold.  Finally,
applying Proposition~\ref{prop: recovery of gen Finler str} generates
the desired Finsler structure.
\end{proof}

\bibliographystyle{amsplain}

\end{document}